\documentclass[11pt]{article}
\usepackage{IEEEtrantools}
\usepackage{mathtools, nccmath}
\usepackage{latexsym}
\usepackage{amsfonts}
\usepackage{amssymb}
\usepackage{psfrag}
\usepackage{graphicx}
\usepackage{epsfig}
\usepackage{amsmath,amsfonts,amsthm}
\usepackage{mathrsfs}
\usepackage{enumerate}
\usepackage{enumitem} 
\usepackage{dsfont}
\usepackage{hyperref}
\usepackage{pst-node}
\usepackage{cleveref}
\usepackage{enumitem}
\usepackage{tikz-cd}
\usepackage{multirow}
\usepackage{tikz}
\usetikzlibrary{arrows,decorations.markings}
\usepackage{bookmark}
\usepackage{hyperref}
\hypersetup{
     colorlinks   = true,
     citecolor    = blue
}

\newcommand\restr[2]{{
  \left.\kern-\nulldelimiterspace 
  #1 
  \littletaller 
  \right|_{#2} 
  }}

\newcommand{\littletaller}{\mathchoice{\vphantom{\big|}}{}{}{}}

\newcommand{\be}{\begin{equation}}
\newcommand{\ee}{\end{equation}}
\newcommand{\bea}{\begin{eqnarray}}
\newcommand{\eea}{\end{eqnarray}}
\newcommand{\bean}{\begin{eqnarray*}}
\newcommand{\eean}{\end{eqnarray*}}
\newcommand{\brray}{\begin{array}}
\newcommand{\erray}{\end{array}}
\newcommand{\biearray}{\begin{IEEEarray}{rCl}}
\newcommand{\eiearray}{\end{IEEEarray}}

\newtheorem{dfn}{Definition}[section]
\newtheorem{thm}[dfn]{Theorem}
\newtheorem{lmma}[dfn]{Lemma}
\newtheorem{ppsn}[dfn]{Proposition}
\newtheorem{crlre}[dfn]{Corollary}
\newtheorem{xmpl}[dfn]{Example}
\newtheorem{rmrk}[dfn]{Remark}

\newcommand{\bdfn}{\begin{dfn}\rm}
\newcommand{\bthm}{\begin{thm}}
\newcommand{\blmma}{\begin{lmma}}
\newcommand{\bppsn}{\begin{ppsn}}
\newcommand{\bcrlre}{\begin{crlre}}
\newcommand{\bxmpl}{\begin{xmpl}}
\newcommand{\brmrk}{\begin{rmrk}\rm}

\newcommand{\edfn}{\end{dfn}}
\newcommand{\ethm}{\end{thm}}
\newcommand{\elmma}{\end{lmma}}
\newcommand{\eppsn}{\end{ppsn}}
\newcommand{\ecrlre}{\end{crlre}}
\newcommand{\exmpl}{\end{xmpl}}
\newcommand{\ermrk}{\end{rmrk}}

\newcommand{\bbc}{\mathbb{C}}

\newcommand{\tr}{\mathrm{tr}}

\def \qed { \mbox{}\hfill
$\Box$\vspace{1ex}}

\title{Sections and Chapters}

\begin{document}
	
	
	\author{\sc{Biplab Pal}}

	\title{Depth 2 inclusions of simple $C^*$-algebras and their weak $C^*$-Hopf algebra symmetries}
	\maketitle
	
	
	\begin{abstract}
Let \( B \subset A \) be a depth $2$ inclusion of simple unital \(C^*\)-algebras with a conditional expectation of index-finite type. We show that the second relative commutant \( B' \cap A_1 \) carries a canonical structure of a weak \(C^*\)-Hopf algebra. Furthermore, we construct an action of this weak \(C^*\)-Hopf algebra on \(A\) for which \(B\) is precisely the fixed-point subalgebra, and we prove that the first basic construction \(A_1\) is isomorphic to the crossed product \( A \rtimes (B' \cap A_1) \). This provides a \(C^*\)-algebraic counterpart of the duality between depth~2 subfactors and weak Hopf algebra symmetry, extending the Ocneanu--Nikshych--Vainerman theory beyond the $II_1$ factor setting.
\end{abstract}

	\bigskip
	
	{\bf AMS Subject Classification No.:} {\large 46}L{\large 37}\,, {\large 47}L{\large 40}\,, {\large 46}L{\large 05}\,, {\large 46}L{\large 35}\,, {\large 46}L{\large 67}\,.
	
	\bigskip
	\hypersetup{linkcolor=blue}
	
	\section{Introduction}
Let $N \subset M$ be an inclusion of type $II_1$ factors with finite Jones index, and consider the associated Jones tower of basic constructions (see \cite{J}):
\[
N \subset M \subset M_1 \subset M_2 \subset \cdots \subset M_n \subset \cdots.
\]
It has been well understood since the early development of subfactor theory that the relative commutants $N' \cap M_k$ and $M' \cap M_{k+1}$ possess remarkably rich algebraic structures, playing a fundamental role in the advancement of the theory. Ocneanu first announced, and it was later proved in \cite{D,L,S}, that for a finite index, depth $2$, irreducible inclusion $N \subset M$ of type $II_1$ factors, the second relative commutant $N' \cap M_1$ admits the structure of a finite-dimensional Kac algebra, moreover, there exists a canonical outer action of $N' \cap M_1$ on $M$ such that the fixed-point subalgebra is $N$, and $M_1$ is isomorphic to the crossed product $M \rtimes (N' \cap M_1)$. The case of irreducible depth $2$ inclusions of infinite index was subsequently investigated in \cite{EN,FI}. In order to generalize the Hopf algebraic framework, B\"ohm and Szlach\'anyi \cite{BS} (and later B\"ohm, Nill, and Szlach\'anyi in \cite{BNS}) introduced the concept of weak Hopf algebras to capture the symmetries arising from certain non-irreducible inclusions of $II_1$ factors. Subsequently, Nikshych and Vainerman \cite{NV} demonstrated that for any finite-index, depth $2$ inclusion $N \subset M$ of $II_1$ factors, the second relative commutant $N' \cap M_1$ naturally acquires the structure of a weak $C^*$-Hopf algebra. Also, this structure gives rise to an action of the weak $C^*$-Hopf algebra $N' \cap M_1$ on $M$, under which $N$ appears as the fixed-point subalgebra and $M_1$ realizes the crossed product $M \rtimes (N' \cap M_1)$. In the case of arbitrary depth $2$ inclusions, possibly of infinite index, Enock and J.-M. Vallin developed a comparable description in terms of pseudo-multiplicative unitaries; see \cite{EV}. Later, Kadison and Nikshych \cite{KN}, obtained analogous results for certain symmetric Markov extensions in a purely algebraic setting.

The framework of inclusions of simple $C^*$-algebras provides a natural extension of the theory of subfactors, encompassing both type $II_1$ and type $III$ cases. Motivated by the recent progress in the classification theory of $C^*$-algebras (\cite{BakshiVedlattice, ER, I, M, R}), the study of such inclusions has attracted considerable attention in recent years. In this article, we consider a unital inclusion of simple $C^*$-algebras $B \subset A$ equipped with a conditional expectation of finite Watatani index, which generalizes the notion of the Jones index. Consider the corresponding tower of basic constructions:
\[
B \subset A \subset A_1 \subset A_2 \subset \cdots \subset A_n \subset \cdots.
\]
In this paper, we extend the results of Nikshych and Vainerman \cite{NV} to a more general setting of inclusions of simple unital $C^*$-algebras. In particular, we establish that the second relative commutant admits the structure of a weak $C^*$-Hopf algebra, thereby providing a unified framework that encompasses and generalizes the case of type $II_1$ factors. Moreover, the underlying \(C^*\)-algebra structure of the resulting weak \(C^*\)-Hopf algebra need not coincide with the usual \(C^*\)-algebra structure of the relative commutant \(B'\cap A_1\) (see Remark~\ref{discussion about difference}). We obtain the following result:
\medskip

\noindent\textbf{Theorem A:} (See Theorem \ref{iiimpfor mainth}) Let $B\subset A$ be an inclusion of simple unital $C^*$-algebras with a conditional expectation of index-finite type. Furthermore, if the inclusion is of depth $2$, then the second relative commutant $B'\cap A_1$ admits a weak $C^*$-Hopf algebra structure.
\medskip

Izumi \cite{I} showed that for any irreducible, depth $2$ inclusion $B \subset A$ of simple unital $C^*$-algebras admitting a conditional expectation of index-finite type, there exists an action of a finite-dimensional Kac algebra $H$ on $A$ such that $B$ is the fixed-point subalgebra. In this paper, we extend Izumi's result to a more general setting of arbitary depth $2$ inclusion of simple unital $C^*$-algebras. While a similar result has been mentioned in \cite{SPreprint}, the proofs provided therein are only outlined, with the expectation that detailed arguments would appear in future. To the best of our knowledge, such a comprehensive exposition has not yet been made available in the literature. Moreover, the approach in \cite{SPreprint} is presented primarily in categorical language. For the sake of completeness and mathematical rigor, we include a full treatment of the result in the present paper. In the setting of simple unital $C^*$-algebras, we obtain the following theorem.

\medskip
\noindent\textbf{Theorem B:} (See Proposition \ref{MAINTHM2} and Theorem \ref{MAINTHM3}) Let $B \subset A$ be an inclusion of simple unital $C^*$-algebras with a conditional expectation of index-finite type. Assume further that the inclusion is of depth $2$, then there exists an action of the weak $C^*$-Hopf algebra $ B' \cap A_1$ on $A$ such that $B$ is the fixed-point subalgebra and $A_1 \cong A \rtimes (B'\cap A_1)$.

As a consequence, for any regular inclusion \( B \subset A \) of simple unital \( C^* \)-algebras (i.e., the unitary normalizer of \( B \) in \( A \) generates \( A \) as a \( C^* \)-algebra) with a conditional expectation of index-finite type, there exists an action of the weak \( C^* \)-Hopf algebra \( B'\cap A_1 \) on $A$ such that $B$ is the corresponding fixed-point subalgebra and $
A_1$ is canonically isomorphic to the corresponding crossed product algebra.

\medskip
The article is structured as follows. \Cref{preliminaries} collects the required preliminaries, including Watatani’s $C^*$-index theory, elements of $C^*$-Fourier theory, and basic notions from the theory of weak $C^*$-Hopf algebras. In \Cref{Duality between relative commutants and coalgebra structures}, we study the duality between the relative commutants $B' \cap A_1$ and $A' \cap A_2$ arising from an inclusion $B \subset A$ of simple unital $C^*$-algebras equipped with a conditional expectation of index-finite type. In \Cref{Weak $C^*$-Hopf algebra structure on relative commutants}, assuming that the inclusion has depth $2$, we prove that the second relative commutant carries a natural weak $C^*$-Hopf algebra structure. Finally, in \Cref{Action and crossed product}, we establish the existence of an action of $B'\cap A_1$ on $A$, for which $B$ appears as the associated fixed-point subalgebra and identify the first basic construction with the corresponding crossed product algebra.

\section{Preliminaries}\label{preliminaries}
In this section, we briefly recall the fundamental notions that will be used throughout the paper. Detailed discussions are omitted, and references included for readers seeking additional details.

\subsection{$C^*$-index theory}
	Let $B\subset A$ be an inclusion of simple unital $C^*$-algebras with a conditional expectation of index-finite type. Watatani has shown \cite{Watataniindex} that such an inclusion admits a unique minimal conditional expectation. The Watatani index of $B\subset A$, denoted by $[A : B]_0$, is defined to be the index of this minimal conditional expectation. Consider the corresponding tower of $C^*$-basic construction:
	\begin{center}
		$B \subset A \subset A_1 \subset A_2 \subset \cdots \subset A_n \subset \cdots$
	\end{center}
	where each inclusion $A_{n-1}\subset A_n$ admits a unique (dual) minimal conditional expectation $E_{n}: A_{n}\to A_{n-1}$ for $n \geq 0$. We adopt the convention $A_{-1}:=B$ and $A_0:=A$ and denote by $e_n$ the Jones projection in $A_n$ for each $n\geq 1$. Throughout, we set $ [A:B]^{-1}_0 = \tau$.
	
	We now recall a few definitions and results that will be used in the subsequent discussion.

	\begin{lmma}[\cite{BakshiVedlattice}]\label{pushdown}
		For any $x_1\in A_1$, there exists a unique element $x_0\in A$ satisfying $x_1e_1=x_0e_1$, where $x_0=\tau^{-1} E_1(x_1e_1)$.
	\end{lmma}
	
	Let $B^\prime\cap A_n:= \{x\in A_n: xb=bx\mbox{ for all }b\in B\}$ denote the relative commutant of $B$ in $A_n$. It is well known that each relative commutant $B'\cap A_n$ is finite dimensional \cite{Watataniindex}. Using the minimal conditional expectations, one can obtain a consistent `Markov type trace' on each $B'\cap A_n$ (\cite{BakshiVedlattice}[Proposition $2.21$]). More concretely, for every $n\geq 0$, the map $\mathrm{tr}_n: B^{\prime}\cap A_n\to\mathbb{C}$ defined by $\mathrm{tr}_n=(E_0\circ E_1\circ\cdots\circ E_n)_{|_{B^{\prime}\cap A_n}}$ is a faithful tracial state on $B^{\prime}\cap A_n$. 
	
	\begin{ppsn}[\cite{BakshiVedlattice}]\label{m1}
		For every $n \geq 0$, the relative commutant $B^{\prime}\cap A_n$ admits a faithful tracial state $\mathrm{tr}_n$ satisfying
		\begin{equation}\nonumber
		\mathrm{tr}_n(xe_n)={\tau}\mathrm{tr}_{n-1}(x) \,\,\text{ for all } x\in B^{\prime}\cap A_{n-1}
		\end{equation}
	and $\restr{\tr_n}{B'\cap A_{n-1}}=\tr_{n-1}$ for all $n\geq 1$. 
	\end{ppsn}
	
	For convenience, we shall omit the subscript and write simply $\tr_n$ by $\tr$. Let $\{\lambda_i:1\leq i\leq n\}\subset A$ denote a quasi-basis for the minimal conditional expectation $E_0:A\to B$. Then, the $\mathrm{tr}$-preserving conditional expectation from $B'\cap A_n$ onto $A'\cap A_n$ is given by, 
	\begin{equation}\label{f2}
	E^{B^{\prime}\cap A_n}_{A^{\prime}\cap A_n}(x)=\tau\sum_{i}\lambda_ix\lambda^*_i\,,\quad \text{ for all }x \in B'\cap A_n.
	\end{equation}
The $\mathrm{tr}$-preserving conditional expectation from $B^{\prime}\cap A_n$ onto $A^{\prime}_1\cap A_n$ is given by
\begin{equation}\label{conditionalexpimp}
E^{B^{\prime}\cap A_n}_{A^{\prime}_1\cap A_n}(x)=\tau\sum_{i,j}\lambda_ie_1\lambda_j\,x\,\lambda^*_je_1\lambda^*_i\,,\quad  \text{ for all } x \in B'\cap A_n.
\end{equation} 
 It is straightforward to verify that
	$E^{B^{\prime}\cap A_1}_{A^{\prime}\cap A_1}(e_1)=\tau.$ Consequently, $\tr(we_1w')=\tau\tr(ww')$ for all $w,w' \in A'\cap A_n$. The following result will be frequently employed in what follows.
	\begin{lmma}[\cite{KajiwaraWatatani}] \label{kajiwara watatani}
     Consider a unital inclusion $B\subset A$ of simple $C^*$-algebras and let $E_0:A \rightarrow B$ be the minimal conditional expectation. Then for every $x \in B' \cap A $ and $a \in A$, the following relation holds, $E_0(xa)=E_0(ax)$.
 \end{lmma} 
 An inclusion \( B \subset A \) of unital \( C^* \)-algebras with a conditional expectation of index-finite type is said to have finite depth (see \cite{JOPT}) if there exists an integer \( n \geq 1 \) satisfying $(B' \cap A_{n-1})\, e_n\, (B' \cap A_{n-1}) = B' \cap A_n.$ The least such integer \( n \) is termed the depth of the inclusion.
 \subsection{$C^*$-Fourier theory}
 Analogous to subfactor theory, in \cite{BakshiVedlattice} the authors have provided a Fourier theory using Watatani’s notions of index and $C^*$-basic construction of certain inclusions of $C^*$-algebras, and it has found significant applications in various contexts (see, for instance, \cite{BakshiVedlattice,BGP,BGS}). The concept of a Fourier transform for paragroups associated with finite depth subfactors was originally introduced by Ocneanu \cite{O}, where it plays a central role in the advancement of subfactor theory. Subsequently, Bisch \cite{B} provided an explicit formula for the Fourier transform on higher relative commutants in the case of extremal subfactors. We begin by recalling the definitions of the Fourier transform and the rotation operators on the relative commutants of an inclusion of simple unital $C^*$-algebras as defined in \cite{BakshiVedlattice}. For each $n\geq 0$, the Fourier transform $\mathcal{F}_n:B^{\prime}\cap A_n\longrightarrow A^{\prime}\cap A_{n+1}$ is given by,
		\begin{center}
			$\mathcal{F}_n(x)=\tau^{-\frac{n+2}{2}}\,E^{B^{\prime}\cap A_{n+1}}_{A^{\prime}\cap A_{n+1}}(xv_{n+1}),\, \quad \text{ for all } x \in B'\cap A_n,$
		\end{center}
where $v_n=e_n\cdots e_1$. Correspondingly, the inverse Fourier transform ${\mathcal{F}}^{-1}_n: A^{\prime}\cap A_{n+1}\longrightarrow B^{\prime}\cap A_n$ is defined as,
		\begin{center}
			$\mathcal{F}^{-1}_n(w)=\tau^{-\frac{n+2}{2}}\,E_{n+1}(wv^*_{n+1}),\, \quad \text{ for all } w\in A'\cap A_{n+1}.$
		\end{center}
The designation `inverse' is appropriate because, for every $n\geq 0$,  $\mathcal{F}_n \circ \mathcal{F}^{-1}_n ={\mathrm{id}}_{A^\prime \cap A_{n+1}} $ and $\mathcal{F}^{-1}_n \circ \mathcal{F}_n ={\mathrm{id}}_{B^\prime \cap A_{n}}$, as shown in (\cite[Theorem 3.5]{BakshiVedlattice}). For every $n \geq 0$, the rotation map ${\rho^+_n}:B^{\prime}\cap A_n\rightarrow B^{\prime}\cap A_n$ is given by
	\begin{equation}\nonumber
	{\rho^+_n}(x)={\big({\mathcal{F}}^{-1}_n\big({\mathcal{F}_n(x)
		}^*\big)\big)}^*, \quad \text{ for all}\quad  x \in B'\cap A_n.
	\end{equation}
In a similar fashion, ${\rho^-_n}: A^\prime \cap A_{n+1} \longrightarrow A^\prime \cap A_{n+1}$ is defined as follows,
	\begin{center}
		${\rho^-_n}(w)=\big(\mathcal{F}_n\big({\mathcal{F}_n^{-1}(w)}^*\big)\big)^* , \quad \text{ for all }w \in A^{\prime}\cap A_{n+1}.$
	\end{center}
Reflection operators, introduced in \cite{BGPS}, yield various Fourier-theoretic inequalities on the higher relative commutants and establish a connection between the Connes–St\o rmer entropy of the canonical shift and the minimal Watatani index. For every $n\geq 0$, the reflection operator ${r^+_{2n+1}}: B^\prime \cap A_{2n+1} \to B^\prime \cap A_{2n+1}$ is given by $r^+_{2n+1}=(\rho^+_{2n+1})^{n+1}.$ Analogously, $r^-_{2n+1}: A^\prime \cap A_{2n+2} \to A^\prime \cap A_{2n+2}$ is defined by $r^-_{2n+1}=(\rho^-_{2n+1})^{n+1}.$ These operators have the following properties:
\begin{thm}[\cite{BGPS}]\label{aboutreflectionoperator}
Let $B\subset A$ be an inclusion of simple unital $C^*$-algebras with a conditional expectation of index-finite type. Then, for every $n\geq 0$, the reflection operators $r^{\pm}_{2n+1}$ are  unital, involutive, $*$-preserving anti-homomorphisms. Moreover each of them preserves the trace.
\end{thm}
Based on the preceding definition, the reflection operators $r^{\pm}_1$ admit the following explicit forms, which will be essential in the subsequent analysis:
$$r^+_1(x)=\tau^{-1}\sum_{i}E_1(e_1 \lambda_ix)e_1\lambda^*_i, \quad \text{ for all } x\in B'\cap A_1,$$
and
$$r^-_1(w)=\tau^{-2}\sum_{i}\lambda_{i}e_1e_2E_2(we_1e_2\lambda^*_i), \quad \text{ for all } w\in A^{\prime}\cap A_2.$$
 \subsection{Weak $C^*$-Hopf Algebra}
We briefly recall some preliminaries related to weak $C^*$-Hopf algebras; for further details, see \cite{BNS,NV,NSW}.
	
 A weak bialgebra is a quintuple \( ({\mathcal{P}}, m, \eta, \Delta, \varepsilon) \) where \( ({\mathcal{P}}, m, \eta) \) is an algebra and \( ({\mathcal{P}}, \Delta, \varepsilon) \) is a coalgebra, satisfying the following compatibility conditions between the two structures:
        \begin{enumerate}
            \item[$(i)$] \( \Delta \) is an algebra homomorphism (not necessarily unital),
            \item[$(ii)$] $\varepsilon(yy'y'') = \varepsilon(y y^{\prime}_{(1)}) \,\varepsilon(y^{\prime}_{(2)} y'') = \varepsilon(y y^{\prime}_{(2)}) \,\varepsilon(y^{\prime}_{(1)} y'')$, for all \( y, y', y''\in {\mathcal{P}} \),
            \item[$(iii)$] The element \( \Delta(1) \in {\mathcal{P}} \otimes {\mathcal{P}} \) satisfies
            \[
            (\mathrm{id} \otimes \Delta)\circ\Delta(1)=(\Delta(1)\otimes 1)\big(1\otimes \Delta(1)\big)=\big(1\otimes \Delta(1)\big)\big(\Delta(1)\otimes 1\big).
            \]
        \end{enumerate}
 A weak Hopf algebra is a weak bialgebra \( ({\mathcal{P}}, m, \eta, \Delta, \varepsilon) \) endowed with a linear map \( S : {\mathcal{P}} \rightarrow {\mathcal{P}} \), referred to as the antipode, which satisfies the following relations for all $y\in {\mathcal{P}}$,
    \begin{enumerate}
        \item[$(i)$] \( y_{(1)} S(y_{(2)}) = \varepsilon(1_{(1)} y) 1_{(2)} \),
        \item[$(ii)$] \( S(y_{(1)}) y_{(2)} = 1_{(1)} \varepsilon(y 1_{(2)}) \),
        \item[$(iii)$] \( S(y_{(1)}) y_{(2)} S(y_{(3)}) = S(y) \).
    \end{enumerate}
A weak \( C^* \)-Hopf algebra is a weak Hopf algebra \( ({\mathcal{P}}, m, \eta, \Delta, \varepsilon) \) for which \( {\mathcal{P}}\) is a finite-dimensional \( C^* \)-algebra and the comultiplication $\Delta$ is \( * \)-preserving, that is, $\Delta(y^*) = {\Delta(y)}^*$, for all $y\in {\mathcal{P}}$. A weak Kac algebra refers to a weak \( C^* \)-Hopf algebra \( ({\mathcal{P}}, m, \eta, \Delta, \varepsilon, S) \) in which the antipode $S$ acts as a $*$-preserving involution, that is, $S^2 = \mathrm{id}$ and $ S(y^*)={S(y)}^*$ for all $y\in {\mathcal{P}}$. We adopt Sweedler’s notation throughout this paper, for all $y\in {\mathcal{P}}$ we write
\[
\Delta(y) = y_{(1)} \otimes y_{(2)}, \quad (\Delta \otimes \mathrm{id})\circ \Delta(y) = y_{(1)} \otimes y_{(2)} \otimes y_{(3)} = (\mathrm{id} \otimes \Delta)\circ \Delta(y).\]
For a weak $C^*$-Hopf algebra \({\mathcal{P}}\), the target and source counital maps, denoted respectively by
\(\varepsilon^t\) and \(\varepsilon^s\) , are defined as
\[
\varepsilon^t(y) = \varepsilon(1_{(1)} y)\,1_{(2)}, 
\quad 
\varepsilon^s(y) = 1_{(1)} \,\varepsilon(y 1_{(2)}),
\quad \text{ for all }  y \in {\mathcal{P}}.
\]
The images of these maps are unital \(C^*\)-subalgebras of \({\mathcal{P}}\), referred to as the Cartan subalgebras:
\[
{\mathcal{P}}_t := \{\, y \in {\mathcal{P}} : \varepsilon^t(y) = y \,\}, 
\quad 
{\mathcal{P}}_s := \{\, y \in {\mathcal{P}} : \varepsilon^s(y) = y \,\}.
\]
We say that \({\mathcal{P}}\) is connected if the inclusion 
\({\mathcal{P}}_t \subset {\mathcal{P}}\) is connected. The dual of a weak $C^*$-Hopf algebra admits a canonical weak $C^*$-Hopf algebra structure. ${\mathcal{P}}$ is connected iff ${\mathcal{P}}^*_t \cap {\mathcal{P}}^*_s= \bbc$ (\cite{N}). We say that \({\mathcal{P}}\) is biconnected if both \({\mathcal{P}}\) and its dual are connected.
\smallskip

The notion of action of a weak $C^*$-Hopf algebra, crossed product and fixed point algebra were introduced in \cite{NSW}. A (left) action of a weak \(C^*\)-Hopf algebra \({\mathcal{P}}\) on a $C^*$-algebra $A$ is a linear map
\[
{\mathcal{P}} \otimes A \;\to\; A, \qquad y \otimes a \;\mapsto\; y \triangleright a,
\]
turning \(A\) into a left \({\mathcal{P}}\)-module such that, for each $y\in {\mathcal{P}}$ the operator $y \otimes a \;\mapsto\; y \triangleright a$ is continuous and satisfies the following relations for all \(y \in {\mathcal{P}}\) and \(a,a' \in A\), 
\begin{enumerate}
\item[$(i)$] \(y \triangleright (aa') = (y_{(1)} \triangleright a)\,(y_{(2)} \triangleright a'),\)
\item[$(ii)$] \((y \triangleright a)^* = {S(y)}^* \triangleright a^*,\)
\item[$(iii)$] \(y \triangleright 1 = \varepsilon^t(y) \triangleright 1,\) and  
      \(y \triangleright 1 = 0 \iff \varepsilon^t(y) = 0,\).
\end{enumerate}
Given such an action, one can construct the crossed product algebra \( A \rtimes {\mathcal{P}} \) as follows. As a vector space, \( A \rtimes {\mathcal{P}} \) is defined to be the relative tensor product \( A \otimes_{{\mathcal{P}}_t} {\mathcal{P}} \). Here, \( {\mathcal{P}} \) is regarded as a left \( {\mathcal{P}}_t \)-module under multiplication, whereas \( A \) becomes a right \( {\mathcal{P}}_t \)-module via $a \cdot u = a (u \triangleright 1)$, for all $a \in A,\, u \in {\mathcal{P}}_t.$ We identify elements according to the relation 
\[
a (u \triangleright 1) \otimes y \;\sim\; a \otimes u y,
\qquad a \in A,\, y \in {\mathcal{P}},\, u \in {\mathcal{P}}_t.
\]
For \( a \in A \) and \( y \in {\mathcal{P}} \), let \( [a \otimes y] \) denote the equivalence class of \( a \otimes y \). The multiplication and involution on \( A \rtimes {\mathcal{P}} \) are then defined by 
\[
[a \otimes y][a' \otimes y'] 
= [\, a (y_{(1)} \triangleright a') \otimes y_{(2)} y' \,],
\qquad
{[a \otimes y]}^* 
= [\, (y_{(1)}^* \triangleright a^*) \otimes y_{(2)}^* \,],
\]
for all \( a, a' \in A \) and \( y, y' \in {\mathcal{P}} \). 
It was shown in \cite{NSW} that the algebra \( A \rtimes {\mathcal{P}} \), defined abstractly as above, can be realized as a \( C^* \)-algebra. The $C^*$-subalgebra $A^{\mathcal{P}}=\{a\in A : y\triangleright a = \varepsilon^t(y) \triangleright a, \text{ for all } y\in {\mathcal{P}}\}$ of $A$, is called a fixed point subalgebra. We say that the action $\triangleright$ is minimal if ${\mathcal{P}}_s= A' \cap (A \rtimes {\mathcal{P}}),$ where \(A\) is regarded as a \(C^*\)-subalgebra of \(A\rtimes\mathcal{P}\) via the embedding \(A\ni a\mapsto[a\otimes1]\in A\rtimes\mathcal{P}\).

\section{Duality between relative commutants and coalgebra structures}\label{Duality between relative commutants and coalgebra structures}
In this section, we establish that for any inclusion $B\subset A$ of simple unital $C^*$-algebras with a conditional expectation of index-finite type, there exists a natural duality between the relative commutants $B'\cap A_1$ and $A'\cap A_2$. This duality generalizes a result previously proved for inclusions of $II_1$ factors in \cite{NV}, where it was obtained only for inclusions of depth $2$. In contrast, we prove the duality for arbitrary inclusions. Also, a key distinction from the type $II_1$ case is that modular conjugation operator is not available in our setting; instead, we exploit the quasi-basis and minimal conditional expectation to circumvent this difficulty. Using this duality, we construct coalgebra structures on the relative commutants arising from arbitrary inclusions of simple unital $C^*$-algebras.

\smallskip
\noindent\textbf{Notation:} Throughout the paper, we consider an inclusion $B\subset A$ of simple unital $C^*$-algebras equipped with a conditional expectation of index-finite type. Denote by $\{\lambda_i:1\leq i\leq n\}\subset A$ a quasi-basis for the minimal conditional expectation $E_0$. We also write
 $$\mathcal{P}=B'\cap A_1 \quad\text{ and } \quad \mathcal{Q}=A'\cap A_2.$$
We begin by establishing the following result, which will play a crucial role in the sequel.

\begin{lmma}\label{non degenerate}
We have $(B'\cap A_2)\,e_2= \mathcal{P}e_2$ and $(B'\cap A_2)\,e_1=\mathcal{Q}e_1$.
\end{lmma}
\begin{prf}
Using \Cref{pushdown}, the first part follows immediately. Now let $x\in B' \cap A_2$. Then
\begin{eqnarray}\nonumber
\tau^{-1}E^{B^{\prime}\cap A_2}_{A^{\prime}\cap A_2}(xe_1)e_1&=&\sum_{i}\lambda_ixe_1\lambda^*_ie_1
=\sum_{i}\lambda_ix E_0(\lambda^*_i)e_1= x e_1.
\end{eqnarray}
Hence $(B'\cap A_2)\,e_1\subset \mathcal{Q}e_1$. Since it is clear that $\mathcal{Q}e_1 \subset (B'\cap A_2)\,e_1$, the result follows.\qed
\end{prf}

The following proposition is motivated by \cite[Proposition 3.2]{NV}.
\begin{ppsn}\label{form}
The following non-degenerate bilinear form establishes a non-degenerate duality between $\mathcal{P}$ and $\mathcal{Q}$:
 $$\big<x,w\big>_{\text{NV}}=d\tau^{-2}\tr(xe_2e_1w),\quad x\in \mathcal{P} \text{ and } w\in \mathcal{Q},$$
where \( d \) denotes the constant number  \( \| {\mathrm{Ind}}_{\text{W}}(E_1|_{A' \cap A_1}) \|_2 \).
\end{ppsn}
\begin{prf} Let $w\in \mathcal{Q}$ be such that $\big<\mathcal{P},w\big>_{\text{NV}}=0$. Then we have
\begin{eqnarray}\nonumber
\tr(ww^*)&=&\tau^{-1}\tr\big(ww^* \,E^{B^{\prime}\cap A_2}_{A^{\prime}\cap A_2}(e_1)\big)\\\nonumber
&=&\tau^{-1}\tr\big(E^{B^{\prime}\cap A_2}_{A^{\prime}\cap A_2}(ww^* e_1)\big)\\\nonumber
&=&\tau^{-1}\tr(ww^* e_1)\\\nonumber
&=&\tau^{-2}\tr(e_1e_2e_1ww^*)\\\nonumber
&=&\tau^{-2}\tr(w^*e_1e_2e_1w).
\end{eqnarray}
By \Cref{non degenerate}, we have $w^*e_1 \in \mathcal{Q}e_1 = (B' \cap A_2)e_1 \subset (B' \cap A_2)$. Let $w^*e_1 = x_2$. Then $x_2e_2 \in (B' \cap A_2)e_2 = \mathcal{P}e_2$. Hence, we have $w^*e_1e_2 = x_2e_2 = x_1e_2$, where $x_1 \in \mathcal{P}$. Thus, we obtain the following:

\begin{eqnarray}\nonumber
\tr(ww^*)
&=&\tau^{-2}\tr(x_1e_2e_1w)=d^{-1}\big<x_1,w\big>_{\text{NV}}=0.
\end{eqnarray}
Hence, we get $w=0$. Now, let $x\in \mathcal{P}$ be such that $\big<x,\mathcal{Q}\big>_{\text{NV}}=0$. Then by \Cref{non degenerate}, $\tr(xe_2e_1(B'\cap A_2))=0.$
It follows that,
\begin{eqnarray}\nonumber
\tr(xx^*)&=& E_0\circ E_1(x^*x)\\\nonumber
&=& \tau^{-1} E_0\circ E_1 \circ E_2(x^*x e_2)\\\nonumber
&=&\tau^{-1} E_0\circ E_1 \circ E_2(x e_2x^*)\\\nonumber
&=&\tau^{-2} E_0\circ E_1 \circ E_2(x e_2e_1e_2x^*)\\\nonumber
&=&0.
\end{eqnarray}
Hence, we obtain $x=0$, and the proof is complete.
\qed
\end{prf}
\begin{dfn}\label{NV}
Using the form $\langle\, , \,\rangle_{\text{NV}}$, we define the comultiplication $\Delta_{\mathcal{Q}}$, the counit $\varepsilon_{\mathcal{Q}}$, and the antipode $S_{\mathcal{Q}}$ on ${\mathcal{Q}}$ as follows:
$$\Delta_{\mathcal{Q}}:{\mathcal{Q}} \rightarrow {\mathcal{Q}} \otimes {\mathcal{Q}}: \quad  \,\big<x_1x_2,w\big>_{\text{NV}}=\big<x_1,w_{(1)}\big>_{\text{NV}} \,\,\big<x_2,w_{(2)}\big>_{\text{NV}}\, ,$$
$$\varepsilon_{\mathcal{Q}}: {\mathcal{Q}} \rightarrow \bbc: \quad \varepsilon_{\mathcal{Q}}(w)=\big<1,w\big>_{\text{NV}}\, ,$$
$$S_{\mathcal{Q}}:\mathcal{Q}\rightarrow \mathcal{Q}: \quad \big<x,S_{\mathcal{Q}}(w)\big>_{\text{NV}}=\overline{\big<x^*,w^*\big>}_{{\text{NV}}}\, , $$
for all $x,x_1,x_2 \in \mathcal{P}$ and $w \in \mathcal{Q}$.
\end{dfn}
Clearly, the triple $(\mathcal{Q} ,\Delta_{\mathcal{Q}},\varepsilon_{\mathcal{Q}}) $ constitutes a coalgebra.

\section{Weak $C^*$-Hopf algebra structure on relative commutants}\label{Weak $C^*$-Hopf algebra structure on relative commutants}

The objective of this section is to establish that, for an inclusion $B\subset A$  of simple unital $C^*$-algebras with a conditional expectation of index-finite type, and moreover, if the inclusion is of depth $2$, then the second relative commutant admits a weak $C^*$-Hopf algebra structure. Furthermore, if  $E_1|_{A'\cap A_1}$ has a scalar Watatani index, then this weak $C^*$-Hopf algebra structure refines to a weak Kac algebra structure. Although our main results concern depth $2$ inclusions, we first record several preliminary observations that hold for an arbitrary inclusion $B\subset A$ of simple unital $C^*$-algebras equipped with a conditional expectation of index-finite type.
 \begin{ppsn}\label{s=r-}
We have the following, $S_{\mathcal{Q}}=r^-_1.$
\end{ppsn}
\begin{prf}
For each $x\in \mathcal{P}$, it holds that
 \begin{eqnarray}\nonumber
 \big<x,r^-_1(w)\big>_{\text{NV}}&=& d\tau^{-2}\tr(xe_2e_1r^-_1(w))\\\nonumber
 &=& d\tau^{-4}\sum_{i}\tr(xe_2e_1\lambda_i e_1 e_2E_2(we_1e_2\lambda^*_i))\\\nonumber
 &=& d\tau^{-3}\sum_{i}\tr(x E_0(\lambda_i) e_2E_2(we_1e_2\lambda^*_i))\\\nonumber
 &=& d\tau^{-3}\tr(x e_2E_2(we_1e_2))\\\nonumber
 &=& d\tau^{-2}\tr(xE_2(we_1e_2))\\\label{for*preserving}
 &=& d\tau^{-2}\tr(xwe_1e_2)\\\nonumber
 &=& \overline{\big<x^*,w^*\big>}_{{\text{NV}}}\,.
 \end{eqnarray}
 Therefore, we obtain the desired conclusion.
 \qed
\end{prf}

From \Cref{aboutreflectionoperator}, we know that $r^-_1$ is an involutive anti-homomorphism. Therefore, we obtain the following result:
\begin{crlre}\label{imp for kac alg 4}
$S_{\mathcal{Q}}$ is an anti-algebra map such that $S_{\mathcal{Q}}^2=\mathrm{id}.$
\end{crlre}

Now we will record the following crucial observations:
\begin{lmma}\label{impr-r+}
 For any $x\in \mathcal{P}$ and $w\in \mathcal{Q}$ we have,
$$\big<x,r^-_1(w)\big>_{\text{NV}}= \big<r^+_1(x),w\big>_{\text{NV}}\,.$$
\end{lmma}
\begin{prf}
We proceed as follows. For any $x\in \mathcal{P}$ and $w\in \mathcal{Q}$ we have,
\begin{eqnarray}\nonumber
\big<r^+_1(x),w\big>_{\text{NV}}&=& d\tau^{-2}\tr(r^+_1(x)e_2e_1w)\\\nonumber
&=& d\tau^{-3}\sum_{i}\tr(e_2e_1wE_1(e_1\lambda_i x)e_1\lambda^*_i)\\\nonumber
&=& d\tau^{-2}\sum_{i}\tr(e_2e_1\lambda_i xwe_1\lambda^*_i)\\\nonumber
&=& d\tau^{-3}\sum_{i}\tr\big(e_1E^{B^{\prime}\cap A_2}_{A^{\prime}\cap A_2}( xwe_1e_2)\big).
\end{eqnarray}
From \cite{BakshiVedlattice} [Corollary 2.24], we know that $E^{B^{\prime}\cap A_2}_{A^{\prime}\cap A_2}(e_1)=\tau.$ Hence we obtain,
\begin{eqnarray}\nonumber
\big<r^+_1(x),w\big>_{\text{NV}}&=& d\tau^{-2}\sum_{i}\tr\big(E^{B^{\prime}\cap A_2}_{A^{\prime}\cap A_2}( xwe_1e_2)\big)= d\tau^{-2}\tr(xwe_1e_2).
\end{eqnarray}
Thus using \Cref{for*preserving}, we are done. \qed
\end{prf}
 \begin{ppsn}\label{Delta*preserving}
 $\Delta_\mathcal{Q}$ is $*$-preserving.
 \end{ppsn}
\begin{prf}
Using \Cref{s=r-} and \Cref{impr-r+}, for any $x_1,x_2 \in \mathcal{Q}$ and $w \in \mathcal{Q}$ we have
\begin{eqnarray}\nonumber
\big<x_1x_2,w^*\big>_{\text{NV}}
&=&\overline{\big<x^*_2x^*_1,r^-_1(w)\big>}_{\text{NV}}\\\nonumber
&=&\overline{\big<r^+_1(x^*_1)r^+_1(x^*_2),w\big>}_{\text{NV}}\\\nonumber
&=&\overline{\big<r^+_1(x^*_1),w_{(1)}\big>}_{\text{NV}} \,\,\overline{\big<r^+_1(x^*_2),w_{(2)}\big>}_{\text{NV}}\\\nonumber
&=&\overline{\big<x^*_1,r^-_1(w_{(1)})\big>}_{\text{NV}} \,\,\overline{\big<x^*_2,r^-_1(w_{(2)})\big>}_{\text{NV}}\\\nonumber
&=&\big<x_1,(w_{(1)})^*\big>_{\text{NV}} \,\,\big<x_2,(w_{(2)})^*\big>_{\text{NV}}\, .
\end{eqnarray}
Consequently, the result follows. \qed
\end{prf}

As an immediate consequence we get the following:
\begin{crlre} \label{imp for kac alg 3}
$S_\mathcal{Q}$ is an anti-coalgebra map.
\end{crlre}
\begin{lmma}\label{property of reflection}
We have, $r^-_1(A'\cap A_1)= A^{\prime}_1\cap A_2.$
\end{lmma}
\begin{prf}
Let $v\in A'\cap A_1$. Then, we obtain the following:
\begin{eqnarray}\nonumber
r^-_1(v)&=&\tau^{-2}\sum_{i}\lambda_{i}e_1e_2E_2(ve_1e_2\lambda^*_i)\\\nonumber
&=&\tau^{-1}\sum_{i}\lambda_{i}e_1e_2ve_1\lambda^*_i\\\nonumber
&=&\tau^{-1}\sum_{i}\lambda_{i}e_1E^{B^{\prime}\cap A_2}_{A^{\prime}\cap A_2}( e_2v)e_1\lambda^*_i\\\nonumber
&=&\sum_{i,j}\lambda_{i} e_1 \lambda_{j} e_2 v\lambda^*_je_1\lambda^*_i\\\nonumber
&=&\tau^{-1} E^{B^{\prime}\cap A_2}_{A^{\prime}_1\cap A_2}(e_2v).\quad \text{ (by \Cref{conditionalexpimp})} 
\end{eqnarray}
Thus we obtain $r^-_1(v)\in A^{\prime}_1\cap A_2$. Now, if $w\in A^{\prime}_1\cap A_2$, then we have
\begin{eqnarray}\nonumber
r^-_1(w)&=&\tau^{-2}\sum_{i}\lambda_{i}e_1e_2E_2(we_1e_2\lambda^*_i)\\\nonumber
&=&\tau^{-2}\sum_{i}\lambda_{i}e_1e_2E_2(e_1\lambda^*_i e_2w) \quad \text{ (by \Cref{kajiwara watatani})}\\\nonumber
&=&\tau^{-1}\sum_{i}\lambda_{i}e_1\lambda^*_i E_2(e_2w)\\\nonumber
&=&\tau^{-1} E_2(e_2w).
\end{eqnarray}
So we have $r^-_1(w)\in A'\cap A_1$. By \Cref{aboutreflectionoperator}, since $r^-_1$ is an involution, we obtain the desired result.
\qed
\end{prf}

 We know that $E_1|_{A'\cap A_1}$ admits a quasi-basis in $A'\cap A_1$; let $\{u_i\}$ denote such a quasi-basis. We now present the following observations, whose proofs are influenced by \cite{KN}.
\begin{lmma}\label{Delta(1)}
$\Delta_\mathcal{Q}(1)=d^{-1} \sum_{i} r^{-}_1(u^*_i) \otimes u_i$.
 
\end{lmma}
\begin{prf}
For all $x,y \in \mathcal{P}$, we have the following:
\begin{eqnarray}\nonumber
d^{-1}\sum_{i}\big<x,r^-_1(u^*_i)\big>_{\text{NV}}\,\big<y,u_i\big>_{\text{NV}}&=&d\tau^{-4}\sum_{i}\tr(u^*_ie_1e_2x)\,\tr(ye_2e_1u_i)\\\nonumber
&=&d\tau^{-4}\sum_{i}\tr(u^*_ie_1e_2x)\,E_1\big(E^{B^{\prime}\cap A_1}_{A^{\prime}\cap A_1}\circ E_2( ye_2e_1u_i)\big)\\\nonumber
&=&d\tau^{-3}\sum_{i}\tr(u^*_ie_1e_2x)\,E_1\big(E^{B^{\prime}\cap A_1}_{A^{\prime}\cap A_1}( ye_1u_i)\big)\\\nonumber
&=&d\tau^{-4}\sum_{i}\tr(u^*_ie_1e_2x)\,E_1\big(E^{B^{\prime}\cap A_1}_{A^{\prime}\cap A_1}( ye_1)u_i\big)\\\nonumber
&=&d\tau^{-3}\tr(E^{B^{\prime}\cap A_1}_{A^{\prime}\cap A_1}( ye_1)e_1e_2x).
\end{eqnarray}
Using \Cref{f2}, we obtain that
$$E^{B^{\prime}\cap A_1}_{A^{\prime}\cap A_1}( ye_1)e_1e_2x=\tau\sum_j \lambda_jye_1\lambda^*_je_1e_2x=\tau ye_1e_2x.$$
So finally we have,
\begin{eqnarray}\nonumber
d^{-1}\sum_{i}\big<x,r^-_1(u^*_i)\big>_{\text{NV}}\, \big<y,u_i\big>_{\text{NV}}&=&d\tau^{-2}\tr(ye_1e_2x)
=d\tau^{-2}\tr(ye_2e_1x)
=\big<xy,1\big>_{\text{NV}}\,.
\end{eqnarray}
Hence we are done.\qed
\end{prf}

\begin{lmma}\label{item (2)Delta(1)}
For any $w\in \mathcal{Q}$ and $v \in A'\cap A_1$, we have the following:
$$\Delta_\mathcal{Q}(wv)=\Delta_\mathcal{Q}(w)(v\otimes 1).$$
\end{lmma}
\begin{ppsn}\label{imp for kac algebra}
 We have the following:
\begin{eqnarray}\nonumber
(\mathrm{id}\otimes \Delta_\mathcal{Q})\circ\Delta_{\mathcal{Q}}(1)&=&\big(\Delta_{\mathcal{Q}}(1)\otimes 1\big)\big(1\otimes \Delta_{\mathcal{Q}}(1)\big)
=\big(1\otimes \Delta_{\mathcal{Q}}(1)\big)\big(\Delta_{\mathcal{Q}}(1)\otimes 1\big).
\end{eqnarray}
\end{ppsn}
\begin{prf}
Using \Cref{Delta(1)}, we have
$$\Delta_{\mathcal{Q}}(1)=d^{-1} \sum_{i} r^{-}_1(u^*_i) \otimes u_i$$
Thus, using \Cref{property of reflection}, we obtain the following:
\begin{eqnarray}\nonumber
\big(\Delta_{\mathcal{Q}}(1)\otimes 1\big)\big(1\otimes \Delta_{\mathcal{Q}}(1)\big)&=& d^{-2} \sum_{i,j}( r^{-}_1(u^*_i) \otimes u_i\otimes 1)(1\otimes r^{-}_1(u^*_j) \otimes u_j)\\\nonumber
&=& d^{-2} \sum_{i,j} r^{-}_1(u^*_i) \otimes u_i r^{-}_1(u^*_j)\otimes u_j\\\label{impinsomeprf1}
&=& d^{-2} \sum_{i,j} r^{-}_1(u^*_i) \otimes  r^{-}_1(u^*_j)u_i\otimes u_j\\\nonumber
&=&\big(1\otimes \Delta_{\mathcal{Q}}(1)\big)\big(\Delta_{\mathcal{Q}}(1)\otimes 1\big).
\end{eqnarray}
Moreover, using \Cref{item (2)Delta(1)}, we deduce that
\begin{eqnarray}\nonumber
(\mathrm{id} \otimes \Delta_{\mathcal{Q}})\circ\Delta_{\mathcal{Q}}(1)&=& d^{-1} \sum_{i} r^{-}_1(u^*_i) \otimes \Delta_{\mathcal{Q}}(1u_i)\\\label{impinsomeprf2}
&=& d^{-2} \sum_{i,j} r^{-}_1(u^*_i) \otimes r^{-}_1(u^*_j)u_i \otimes u_j.
\end{eqnarray}
Finally, from \Cref{impinsomeprf1} and \Cref{impinsomeprf2}, we obtain the desired result. This completes the proof.
\qed
\end{prf}

We now present the following observations, whose straightforward proofs are similar to those in \cite{NV} and are thus omitted.
\begin{lmma}\label{same as NV}
For all $x\in {\mathcal{P}}$ and $w, w' \in {\mathcal{Q}}$, we have the following:
\begin{enumerate}
\item[$(i)$] $E_2(wxe_2)=E_2(e_2x\,r^{-}_1(w)).$
\item[$(ii)$] $\big<x,ww'\big>_{\text{NV}}=\tau^{-1}\big<E_2(w'xe_2),w\big>_{\text{NV}}\,$.
\item[$(iii)$]  $\varepsilon_{\mathcal{Q}}^t(w):=\varepsilon_{\mathcal Q}(1_{(1)}w)1_{(2)}=\tau^{-1}E_2(we_2)$.
\item[$(iv)$] $w\,\varepsilon_{\mathcal{Q}}^t(w')=\varepsilon_{\mathcal{Q}}(w_{(1)}w')\,w_{(2)}$.
\end{enumerate}

\end{lmma}

\begin{ppsn}\label{imp for kac alg 2}
For all $w,w',w''\in {\mathcal{Q}}$, the following holds:
$$\varepsilon_{\mathcal{Q}}(ww'w'')=\varepsilon_{\mathcal{Q}}(ww^{\prime}_{(1)})\varepsilon_{\mathcal{Q}}(w^{\prime}_{(2)}w'')=\varepsilon_{\mathcal{Q}}(ww^{\prime}_{(2)})\varepsilon_{\mathcal{Q}}(w^{\prime}_{(1)}w'').$$
\end{ppsn}
\begin{prf}
We have,
\begin{eqnarray}\nonumber
\varepsilon_{\mathcal{Q}}(ww^{\prime}_{(2)})\varepsilon_{\mathcal{Q}}(w^{\prime}_{(1)}w'')&=&\varepsilon_{\mathcal{Q}}\big(w\,\varepsilon_{\mathcal{Q}}(w^{\prime}_{(1)}w''\big) w^{\prime}_{(2)})\\\nonumber
&=&\varepsilon_{\mathcal{Q}}\big(ww^{\prime} \varepsilon_{\mathcal{Q}}^t(w'')\big) \quad \text{(by item ($iv$) of \Cref{same as NV})} \\\nonumber
&=& d \tau^{-2}\tr\big(e_2e_1ww^{\prime} \varepsilon_{\mathcal{Q}}^t(w'')\big)\\\nonumber
&=& d \tau^{-3}\tr\big(E_2(w''e_2)e_2e_1ww^{\prime}\big) \quad \text{(by item ($iii$) of \Cref{same as NV})}\\\nonumber
&=& d \tau^{-2}\tr(w''e_2e_1ww^{\prime}) \quad \text{(by \Cref{pushdown})}\\\nonumber
&=&\varepsilon_{\mathcal{Q}}(ww'w'').
\end{eqnarray}
By \Cref{Delta*preserving}, we know that $\Delta_{\mathcal{Q}}$ is $*$-preserving.  Moreover, since $({\mathcal{Q}} ,\Delta_{{\mathcal{Q}}},\varepsilon_{{\mathcal{Q}}}) $ is a coalgebra, it follows that $\varepsilon_{\mathcal{Q}}(x) = \overline{\varepsilon_{\mathcal{Q}}(x^*)}$. Therefore, we have
\begin{eqnarray}\nonumber
\varepsilon_{\mathcal{Q}}(ww'w'')&=&\overline{\varepsilon_{\mathcal{Q}}\big((w'')^*(w')^*w^*)\big)} \\\nonumber
&=&\overline{\varepsilon_{\mathcal{Q}}\big((w'')^*(w^{\prime})^*_{(2)}\big)\varepsilon_{\mathcal{Q}}\big((w^{\prime})^*_{(1)}w^*\big)} \\\nonumber
&=&\varepsilon_{\mathcal{Q}}(ww^{\prime}_{(1)})\varepsilon_{\mathcal{Q}}(w^{\prime}_{(2)}w'').
\end{eqnarray}
Hence we are done.
\qed
\end{prf}

We now show that $\Delta_{\mathcal{Q}}$ is an algebra homomorphism. To this end, we require the following:
\begin{lmma}\label{imp form Dleta homo}
For all $w\in {\mathcal{Q}}$ and $v \in A' \cap A_1$, we have the following:
\begin{enumerate}
\item[$(i)$] $r^-_1(e_2)=e_2$.
\item[$(ii)$] $ve_2=r^-_1(v)e_2$.
\item[$(iii)$] $\tau^{-1}E_2(e_2w)=\varepsilon_{\mathcal{Q}}(w1_{(1)})1_{(2)}$.
\item[$(iv)$] $\Delta_{\mathcal{Q}}(w)(1\otimes v)= \Delta_{\mathcal{Q}}(w)(r^-_1(v)\otimes 1)$.
\item[$(v)$] $\varepsilon_{\mathcal{Q}}\big(S_{\mathcal{Q}}(w)\big)=\varepsilon_{\mathcal{Q}}(w)$.
\end{enumerate}
\end{lmma}
\begin{prf}
\begin{enumerate}[leftmargin=*]
\item[($i$)] The following holds:
\begin{eqnarray}\nonumber
r^-_1(e_2)&=&\tau^{-2}\sum_{i}\lambda_{i}e_1e_2E_2(e_2e_1e_2\lambda^*_i)=\tau^{-1}\sum_{i}\lambda_{i}e_1e_2E_2(e_2\lambda^*_i)=\sum_{i}\lambda_{i}e_1e_2\lambda^*_i=e_2.
\end{eqnarray}
\item[($ii$)] Using item ($i$) and \Cref{for*preserving}, for any $x \in {\mathcal{P}}$, we have
\begin{eqnarray}\nonumber
\big<x,r^-_1(v)e_2\big>_{\text{NV}}&=&\big<x,r^-_1(e_2v)\big>_{\text{NV}}= d\tau^{-2}\tr(e_2ve_1e_2x)= d\tau^{-2}\tr(e_2E_1(ve_1)x).
\end{eqnarray}
From \Cref{kajiwara watatani}, we have $E_1(ve_1)=E_1(e_1v)$. Therefore, we obtain
\begin{eqnarray}\nonumber
\big<x,r^-_1(v)e_2\big>_{\text{NV}}&=& d\tau^{-2}\tr(e_2E_1(e_1v)x)=d\tau^{-2}\tr(e_2e_1ve_2x)=\big<x,ve_2\big>_{\text{NV}}\,.
\end{eqnarray}
\item[($iii$)] From \Cref{Delta(1)}, it is easy to observe that, for any $w\in {\mathcal{Q}}$, $\varepsilon_{\mathcal{Q}}(w1_{(1)})1_{(2)} \in A'\cap A_1$. Moreover, it is clear that $E_2(e_2w) \in A'\cap A_1$. So it is enough to show that,
$$\tau^{-1}\tr\big(E_2(e_2w)v\big)=\tr\big(\varepsilon_{\mathcal{Q}}(w1_{(1)})1_{(2)}v\big), \quad \text{for all } v \in A'\cap A_1.$$
We now establish the assertion.
\begin{eqnarray}\nonumber
\tr\big(\varepsilon_	{\mathcal{Q}}(w1_{(1)})1_{(2)}v\big)&=&d\tau^{-2}\tr(e_2e_1w1_{(1)})\tr(1_{(2)}v)\\\nonumber
&=&\tau^{-1}\sum_{i}\tr\big(e_2wr^-_{1}(u^*_i)\big)\tr(u_iv)\\\nonumber
&=&\tau^{-1}\sum_{i}\tr\big(r^-_{1}(u^*_i) e_2w\big)\tr(u_iv)\\\nonumber
&=&\tau^{-1}\sum_{i}\tr(u^*_i e_2 w)\tr(vu_i)\quad \text{(by item ($ii$))}\\\nonumber
&=&\tau^{-1}\tr(ve_2 w)\\\nonumber
&=&\tau^{-1}\tr\big(E_2(e_2w)v\big).
\end{eqnarray}
\item[($iv$)] It is straightforward to verify that,
$$ \big<x,w_{(1)}\,r^-_1(v)\big>_{\text{NV}}\,\big<y,w_{(2)}\big>_{\text{NV}}=\big<x,w_{(1)}\big>_{\text{NV}}\,\big<y,w_{(2)}\,v\big>_{\text{NV}}\,.
$$
\item[($v$)] We have,
 \begin{eqnarray}\nonumber
\varepsilon_{\mathcal{Q}}\big(S_{\mathcal{Q}}(w)\big)&=&\big<1,r^-_1(w)\big>_{\text{NV}}\\\nonumber 
&=& d\tau^{-2}\tr(e_2e_1r^-_1(w))\\\nonumber
&=& d\tau^{-1}\tr(e_2r^-_1(w))\\\nonumber
&=& d\tau^{-1}\tr(e_2w)\\\nonumber
&=& d\tau^{-2}\tr(e_2e_1w)\\\nonumber
&=&\varepsilon_{\mathcal{Q}}(w).
 \end{eqnarray}
\qed
\end{enumerate}
\end{prf}
\begin{lmma}\label{technical1}
Let $\{u_i\}$ be a quasi-basis for $E_1|_{A'\cap A_1}$. Then we have the following:
$$\sum_{i} u^*_i u_i=\sum_{i} u_iu^*_i.$$

\end{lmma}
\begin{prf}
By applying \Cref{same as NV} item ($ii$), for any $x \in {\mathcal{P}}$ we obtain,
\begin{eqnarray}\nonumber
\sum_{i}\big<x,u^*_i u_i\big>_{\text{NV}}&=&\tau^{-1}\sum_i\big<E_2(u^*_i u_ixe_2),1\big>_{\text{NV}}\\\nonumber
&=& \sum_{i}\big<u^*_i u_ix,1\big>_{\text{NV}}\\\nonumber
&=&\sum_{i} \big<u^*_i ,1_{(1)}\big>_{\text{NV}}\,\big< u_ix,1_{(2)}\big>_{\text{NV}}\\\nonumber
&=&d\tau^{-1}\sum_{i} E_1\big(u^*_iE_2(e_21_{(1)})\big)\,\big< u_ix,1_{(2)}\big>_{\text{NV}}.
\end{eqnarray}
Since $E_2(e_21_{(1)}) \in A'\cap A_1$, and $\{u_i\}$ is a quasi-basis of $E_1|_{A'\cap A_1}$. Then we have the following: $$\sum_i u_iE_1\big(u^*_iE_2(e_21_{(1)})\big)=\sum_iE_1\big(E_2(e_21_{(1)})u_i\big)u^*_i=E_2(e_21_{(1)}).$$
Thus we finally have,
\begin{eqnarray}\nonumber
\sum_{i}\big<x,u^*_i u_i\big>_{\text{NV}}&=&d\tau^{-1}\sum_{i} E_1\big(E_2(e_21_{(1)})u_i\big)\,\big<  u^*_ix,1_{(2)}\big>_{\text{NV}}\\\nonumber
&=&d\tau^{-2} \sum_{i}\tr( u_ie_2e_11_{(1)})\,\big<  u^*_ix,1_{(2)}\big>_{\text{NV}}\\\nonumber
&=&\sum_{i}\big<u_iu^*_ix ,1\big>_{\text{NV}}\\\nonumber
&=& \sum_{i}\big<x, u_iu^*_i\big>_{\text{NV}}\,.
\end{eqnarray}
Hence the result follows.
\qed
\end{prf}\\
From now on, we denote $k  ={\mathrm{Ind}}_{{W}}(E_1|_{A'\cap A_1}).$ By applying \Cref{technical1}, we obtain the following proposition.
\begin{ppsn}\label{technical2}
For all $w\in {\mathcal{Q}}$, the following relation holds:
$$\Delta_{\mathcal{Q}}(w)\Delta_{\mathcal{Q}}(1)=d^{-1}\Delta_{\mathcal{Q}}(w)(1\otimes k).$$

\end{ppsn}
\begin{prf}
Using \Cref{Delta(1)} and item $(iv)$ of \Cref{imp form Dleta homo}, we have the following:
\begin{eqnarray}\nonumber
\Delta_{\mathcal{Q}}(w)\Delta_{\mathcal{Q}}(1)&=&d^{-1} \sum_{i}\Delta_{\mathcal{Q}}(w) (r^{-}_1(u^*_i) \otimes u_i)\\\nonumber
&=&d^{-1} \sum_{i}\Delta_{\mathcal{Q}}(w) (1 \otimes u^*_i u_i)\\\nonumber
&=&d^{-1} \sum_{i}(w_{(1)}\otimes w_{(2)}u^*_i u_i).
\end{eqnarray}
Thus from \Cref{technical1}, we get:
 $$\Delta_{\mathcal{Q}}(w)\Delta_{\mathcal{Q}}(1)=d^{-1}\Delta_{\mathcal{Q}}(w)(1\otimes k).$$
\qed
\end{prf}
\begin{lmma}\label{technical3}
For any $x \in {\mathcal{P}}$, $w\in {\mathcal{Q}}$, the following hold:
\begin{enumerate}
\item[$(i)$] $\tau^{-1}E^{B^{\prime}\cap A_2}_{A^{\prime}\cap A_2}(e_1w\,x)=\big<x,w_{(1)}\big>_{\text{NV}}\,w_{(2)}.$
\item[$(ii)$] $\tau^{-1}w_{(2)}E_2(e_2w_{(1)})=d^{-1} wk$.
\end{enumerate}
\end{lmma}
\begin{prf}
 The proof of ($i$) is straightforward. Using item ($iii$) of \Cref{imp form Dleta homo} and \Cref{technical2}, we obtain the following: 
 \begin{eqnarray}\nonumber
 \tau^{-1}w_{(2)}E_2(e_2w_{(1)})&=&\varepsilon_{\mathcal{Q}}(w_{(1)}1_{(1)})w_{(2)} 1_{(2)}\\\nonumber
 &=&\big(\varepsilon_{\mathcal{Q}} \otimes \mathrm{id}\big)\circ\big( \Delta_{\mathcal{Q}}(w)\Delta_{\mathcal{Q}}(1))\\\nonumber
&=&d^{-1} \varepsilon_{\mathcal{Q}}(w_{(1)})w_{(2)}k\\\nonumber
&=&d^{-1} wk.
 \end{eqnarray}
 \qed
\end{prf}

\begin{ppsn}\label{AB=C}
Let $B\subset A$ be an inclusion of simple unital $C^*$-algebras with a conditional expectation of index-finite type. Furthermore, if the inclusion is of depth $2$, then we have, $$B'\cap A_2 = {\mathcal{P}}{\mathcal{Q}}={\mathcal{Q}}{\mathcal{P}}.$$
\end{ppsn}
\begin{prf}
Since the inclusion $B\subset A$ is of depth $2$, it follows from \cite[Proposition 3.11]{JOPT} that $A \subset A_1$ is also of depth 2. Thus we obtain:
$$A'\cap A_3={\mathcal{Q}}e_3{\mathcal{Q}}.$$
Therefore, we obtain $1=\sum_{j}\mu_je_3 \mu^{\prime}_j$, for some $\mu_j, \mu^{\prime}_j  \in {\mathcal{Q}}$. Consequently, for $w \in A_2$, we have:
\begin{eqnarray}\nonumber
we_3&=&\sum_{j}\mu_je_3 \mu^{\prime}_jwe_3\\\nonumber
&=&\sum_{j}\mu_jE_2( \mu^{\prime}_jw)e_3.
\end{eqnarray}
Thus we have $w=\sum_{j}\mu_j E_2( \mu^{\prime}_j w)$. In particular, for any $w \in B'\cap A_2$, we obtain $w \in {\mathcal{Q}}{\mathcal{P}}$. Thus we get $B'\cap A_2 \subset {\mathcal{Q}}{\mathcal{P}}.$ Also for any $w\in   A_2$ we have $w^*=\sum_{j} E_2\big( w^*(\mu^{\prime}_j)^*\big)\mu^*_j.$ Hence, every $w\in B'\cap A_2$ can be expressed as $w=\sum_{j} E_2\big(w (\mu^{\prime}_j)^*\big)\mu^*_j.$ Therefore, we obtain $B'\cap A_2\subset {\mathcal{P}}{\mathcal{Q}}.$ The remaining parts of the inclusions are straightforward, and this concludes the proof.
\qed
\end{prf}

\begin{ppsn}\label{very imp ppsn}
Let $B\subset A$ be an inclusion of simple unital $C^*$-algebras with a conditional expectation of index-finite type. Furthermore, if the inclusion is of depth $2$, then for any $w\in {\mathcal{Q}}$, we have:
$$e_1w=d\tau^{-1} w_{(2)}k^{-1} E_2(e_2e_1w_{(1)}).$$
\end{ppsn}
\begin{prf}
For the sake of completeness, we briefly sketch the proof; similar arguments appear in \cite{KN}. It is straightforward to observe that any $w'\in B'\cap A_2$ satisfying $E^{B^{\prime}\cap A_2}_{A^{\prime}\cap A_2}(w'y)=0$ for all $y\in {\mathcal{P}}$ satisfies $w'=0$. Indeed, $E^{B^{\prime}\cap A_2}_{A^{\prime}\cap A_2}(w'y)=0$ for any $y\in {\mathcal{P}}$ implies $\tr(yww')=0$ for all $y\in {\mathcal{P}}$ and $w\in {\mathcal{Q}}$. Consequently, by applying \Cref{AB=C}, it follows immediately. Hence, it suffices to show that 
\begin{eqnarray}\nonumber
E^{B^{\prime}\cap A_2}_{A^{\prime}\cap A_2}(e_1wy)= d\tau^{-1}E^{B^{\prime}\cap A_2}_{A^{\prime}\cap A_2}( w_{(2)} k^{-1} E_2(e_2e_1w_{(1)})y) \quad \text{for all } y \in {\mathcal{P}}.
\end{eqnarray}
Using item $(i)$ of \Cref{technical3}, we have
\begin{eqnarray}\label{imp}
E^{B^{\prime}\cap A_2}_{A^{\prime}\cap A_2}(e_1w\,y)= \tau\big<y,w_{(1)}\big>_{\text{NV}}\,w_{(2)}.
\end{eqnarray}
Furthermore using \Cref{technical3}, we obtain the following:
\begin{eqnarray}\nonumber
d\tau^{-1}E^{B^{\prime}\cap A_2}_{A^{\prime}\cap A_2}(w_{(2)} k^{-1} E_2(e_2e_1w_{(1)})y)&=&d\tau^{-1} w_{(2)}k^{-1} E^{B^{\prime}\cap A_2}_{A^{\prime}\cap A_2}\big(E_2(e_2e_1w_{(1)}y)\big)\\\nonumber
&=&d\tau^{-1} w_{(2)}k^{-1} E_2\big(e_2E^{B^{\prime}\cap A_2}_{A^{\prime}\cap A_2}(e_1w_{(1)}y)\big)\\\nonumber
&=& d\big<y ,w_{(1)}\big>_{\text{NV}}\, w_{(3)}k^{-1} E_2(e_2 w_{(2)})\\\nonumber
&=& d\big<y ,w_{(1)}\big>_{\text{NV}}\, w_{(3)} E_2(e_2 w_{(2)})k^{-1}\\\label{imp2}
&=&\tau \big<y ,w_{(1)}\big>_{\text{NV}}\,  w_{(2)}. \\\nonumber
\end{eqnarray}
Thus, the result follows from \Cref{imp} and \Cref{imp2}.
\qed
\end{prf}

We now state the following results, whose straightforward proofs are as in Corollary~4.12, Proposition~4.13, Proposition~4.14, and Proposition~4.15 of \cite{NV}, and are therefore omitted.
\begin{crlre}\label{impforaction}
For all $w\in {\mathcal{Q}}$ and $x,y \in A_1$, the following relations hold:
\begin{enumerate}
\item[$(i)$] $xw= d\tau^{-1}  w_{(2)}k^{-1}E_2(e_2xw_{(1)}).$
\item[$(ii)$] $E_2(wxye_2)=d\tau^{-1}E_2(w_{(1)}xe_2)k^{-1}E_2(w_{(2)}ye_2).$
\end{enumerate}
\end{crlre}

\begin{ppsn}\label{delta homo}
For all $w,w' \in {\mathcal{Q}}$, we obtain the followings:
\begin{enumerate}
\item[$(i)$] $\Delta_{\mathcal{Q}}(ww')=\Delta_{\mathcal{Q}}(w)(1 \otimes dk^{-1})\Delta_{\mathcal{Q}}(w')$.
\item[$(ii)$] $d\,w_{(1)}\,r^-_1(w_{(2)}\,k^{-1})=\varepsilon_{\mathcal{Q}}(1_{(1)}w)1_{(2)}.$
\end{enumerate}

\end{ppsn}
\begin{thm}\label{Weak Kac}
Let $B\subset A$ be an inclusion of simple unital $C^*$-algebras with a conditional expectation of index-finite type. Furthermore, if the inclusion is of depth $2$ and $E_1|_{A'\cap A_1}$ has a scalar Watatani index, then $( {\mathcal{Q}}, \Delta_{\mathcal{Q}}, \varepsilon_{\mathcal{Q}}, S_{\mathcal{Q}} )$ is a weak Kac algebra.
\end{thm}
\begin{prf}
Suppose that the restriction $E_1|_{A' \cap A_1}$ has a scalar Watatani index. Then it follows that $d=k$. By invoking \Cref{delta homo} item $(i)$, we see that $\Delta_{\mathcal{Q}}$ is an algebra homomorphism. From \Cref{imp for kac alg 2} and \Cref{imp for kac algebra}, it follows that $({\mathcal{Q}},\, \Delta_{\mathcal{Q}},\, \varepsilon_{\mathcal{Q}})$ is a weak bialgebra. Next, by \Cref{delta homo} item $(ii)$, for any $w \in {\mathcal{Q}}$ we have  
\[
w_{(1)} S_{\mathcal{Q}}(w_{(2)}) 
= \varepsilon_{\mathcal{Q}}(1_{(1)} w)\, 1_{(2)},
\]  
since by \Cref{s=r-}, we have $S_{\mathcal{Q}} = r^-_1$. From \Cref{imp for kac alg 3} and \Cref{imp for kac alg 4}, we know that $r^-_1$ is both an anti-algebra and an anti-coalgebra map, and that it is involutive. It therefore follows that, for any $w \in {\mathcal{Q}}$,
\begin{eqnarray}\nonumber
S_{\mathcal{Q}}(w_{(1)})\, w_{(2)} 
&=& r^-_1\Big(r^-_1(w_{(2)})\, r^-_1\big(r^-_1(w_{(1)})\big)\Big) \\\nonumber
&=& r^-_1\Big(\varepsilon_{\mathcal{Q}}\big(1_{(1)} r^-_1(w)\big)\, 1_{(2)}\Big) \\\nonumber
&=& \varepsilon_{\mathcal{Q}}\big(w\, r^-_1(1_{(1)})\big)\, r^-_1(1_{(2)}) \quad \text{(by \Cref{imp form Dleta homo} item $(v)$)}\\\nonumber
&=& d^{-1}\,\sum_i \varepsilon_{\mathcal{Q}}(w\, u_i^*)\, r^-_1(u_i).\quad \text{(by \Cref{Delta(1)})}
\end{eqnarray}
For any $x\in {\mathcal{P}}$, we have the following:
\begin{eqnarray}\nonumber
\big<x,\,\sum_i \varepsilon_{\mathcal{Q}}(w\, u_i^*)\, r^-_1(u_i)\big>_{\text{NV}}&=&\sum_i \big<1,\, w u_i^*\big>_{\text{NV}}\,\big<x,\, r^-_1(u_i)\big>_{\text{NV}}\\\nonumber
&=&d\tau^{-1}\sum_i\tr(e_2wu^*_i)\big<x,\, r^-_1(u_i)\big>_{\text{NV}}\\\nonumber
&=&d\tau^{-1}\sum_iE_1\big( u^*_iE_2(e_2w)\big)\,\big<x,\, r^-_1(u_i)\big>_{\text{NV}}\\\nonumber
&=&d\tau^{-1}\sum_i\big<x,\, r^-_1(E_1\big( E_2(e_2w)u_i\big)u^*_i)\big>_{\text{NV}}\\\nonumber
&=&\big<x,\,\sum_i \varepsilon_{\mathcal{Q}}(w u_i)\, r^-_1(u_i^*)\big>_{\text{NV}}\,.
\end{eqnarray}
Consequently, we obtain
\begin{eqnarray}\nonumber
S_{\mathcal{Q}}(w_{(1)})\, w_{(2)} 
&=& d^{-1}\, \sum_i\varepsilon_{\mathcal{Q}}(w\, u_i)\, r^-_1(u_i^*)
=1_{(1)}\, \varepsilon_{\mathcal{Q}}(w\, 1_{(2)}).
\end{eqnarray}
Furthermore, since the map $S_{\mathcal{Q}}$ is both an anti-algebra and an anti-coalgebra map, \cite[Lemma~7.4]{Ni} yields
\[
S_{\mathcal{Q}}(w_{(1)})\, w_{(2)}\, S_{\mathcal{Q}}(w_{(3)}) 
= S_{\mathcal{Q}}(w).
\]
Consequently, $( {\mathcal{Q}}, \Delta_{\mathcal{Q}}, \varepsilon_{\mathcal{Q}}, S_{\mathcal{Q}} )$ forms a weak Hopf algebra. By \Cref{Delta*preserving}, this structure is in fact a weak $C^*$-Hopf algebra. Finally, since $S_{\mathcal{Q}} = r^-_1$ is $*$-preserving and satisfies $S^2_{\mathcal{Q}} = \mathrm{id}$, it follows that $( {\mathcal{Q}}, \Delta_{\mathcal{Q}}, \varepsilon_{\mathcal{Q}}, S_{\mathcal{Q}} )$ is a weak Kac algebra.
\qed
\end{prf}

\bigskip

If the Watatani index of $E_1|_{A' \cap A_1}$ is not scalar, then $( {\mathcal{Q}}, \Delta_{\mathcal{Q}}, \varepsilon_{\mathcal{Q}}, S_{\mathcal{Q}} )$ fails to form a weak Kac algebra; for instance, $\Delta_{\mathcal{Q}}$ is no longer a homomorphism. Nevertheless, by introducing a suitable deformation of the structure, along the same lines as in the case of type $II_1$ factors (\cite[Definition~5.1]{NV}), one can obtain a weak $C^*$-Hopf algebra.  

\begin{dfn}\label{dfnnewimp}
We now define the following operations: 
\begin{eqnarray}\nonumber
\text{involution} \quad \dagger:{\mathcal{Q}} \rightarrow {\mathcal{Q}}: \quad w^{\dagger}={r^-_1(k)}^{-1}\,w^*\,r^-_1(k),\\\nonumber
\quad{\tilde{\Delta}}_{\mathcal{Q}} :{\mathcal{Q}} \rightarrow {\mathcal{Q}} \otimes {\mathcal{Q}}:\,\,w_{(\tilde{1})}\otimes w_{(\tilde{2})}
    :=\tilde{\Delta}_{\mathcal{Q}}(w)=d(1\otimes k^{-1})\,\Delta_{\mathcal{Q}}(w),\\\nonumber
\tilde{\varepsilon}_{\mathcal{Q}}:{\mathcal{Q}} \rightarrow \bbc:\quad \tilde{\varepsilon}_{\mathcal{Q}}(w)=d^{-1}\,\varepsilon_{\mathcal{Q}}(kw),\\\nonumber
\tilde{S}_{\mathcal{Q}}:{\mathcal{Q}} \rightarrow {\mathcal{Q}}:\quad \tilde{S}_{\mathcal{Q}}(w)=S_{\mathcal{Q}}(kwk^{-1}).
\end{eqnarray}
\end{dfn}
It is straightforward to verify that the operation $\dagger
$ equips ${\mathcal{Q}}$ with a $C^*$-algebra structure. We next show that, equipped with the above structures, \(\mathcal{Q}\) becomes a weak \(C^*\)-Hopf algebra. The following results from \cite[Lemmas~4.12 and~4.22]{BGoP} will be used repeatedly in the sequel:
\begin{equation}\label{revision1}
r^+_1(B'\cap A)=A'\cap A_1.
\end{equation}
and
\begin{equation}\label{revision2}
e_1r^+_1(u')=e_1u',\qquad \text{for every }u'\in B'\cap A.
\end{equation}
First, we establish the following auxiliary result.
\begin{lmma}\label{revision33}
For any $w \in {\mathcal{Q}}$, we have
$
\tilde{\varepsilon}_{\mathcal{Q}}\!\big(
\tilde{S}_{\mathcal{Q}}(w)\big)
= \tilde{\varepsilon}_{\mathcal{Q}}(w).$ 
\end{lmma}
\begin{prf}
We have the following:
\begin{eqnarray}\nonumber
\tilde{\varepsilon}_{\mathcal{Q}}\!\big(
\tilde{S}_{\mathcal{Q}}(w)\big)
&=&d^{-1}\,\varepsilon_{\mathcal{Q}}(k\,\tilde{S}_{\mathcal{Q}}(w))\\\nonumber
&=& \tau^{-2}\tr(e_2e_1k\,r^-_1(kwk^{-1})).
\end{eqnarray}
By \Cref{revision1}, there exists \(u'\in B'\cap A\) such that
\(k=r_1^+(u')\). Applying \Cref{revision2}, we obtain
\begin{eqnarray}\nonumber
\tilde{\varepsilon}_{\mathcal{Q}}\!\big(
\tilde{S}_{\mathcal{Q}}(w)\big)&=& \tau^{-2}\tr(e_2e_1\,r^+_1(u')\,r^-_1(kwk^{-1}))\\\nonumber
&=& \tau^{-2}\tr(e_2e_1\,r^-_1(kwk^{-1})u')\\\nonumber
&=& \tau^{-2}\tr(e_2u'e_1\,r^-_1(kwk^{-1}))\\\nonumber
&=& \tau^{-2}\tr(e_2ke_1\,r^-_1(kwk^{-1}))\\\nonumber
&=& \tau^{-2}\tr(e_2r_1^-(k)\,e_1\,r^-_1(kwk^{-1})). \quad \text{ (by \Cref{imp form Dleta homo} item $(ii)$)}
\end{eqnarray}
It follows from \Cref{property of reflection} that
\(r_1^-(k)\in A_1'\cap A_2\). Hence,
\begin{eqnarray}\nonumber
\tilde{\varepsilon}_{\mathcal{Q}}\!\big(
\tilde{S}_{\mathcal{Q}}(w)\big)
&=&\tau^{-2}\tr(e_2e_1\,r^-_1(kw))\\\nonumber
&=& d^{-1}\,\varepsilon_{\mathcal{Q}}({S}_{\mathcal{Q}}(kw))\\\nonumber
&=& d^{-1}\,\varepsilon_{\mathcal{Q}}(kw)\quad \text{ (by \Cref{imp form Dleta homo} item $(v)$)}\\\nonumber
&=& \tilde{\varepsilon}_{\mathcal{Q}}(w).
\end{eqnarray}
\qed
\end{prf}
\begin{thm}\label{MAINTHM1}
Let $B\subset A$ be an inclusion of simple unital $C^*$-algebras with a conditional expectation of index-finite type. Furthermore, if the inclusion is of depth $2$, then 
$({\mathcal{Q}}, \tilde{\Delta}_{\mathcal{Q}}, \tilde{\varepsilon}_{\mathcal{Q}}, \tilde{S}_
{\mathcal{Q}})$ forms a weak $C^*$-Hopf algebra.
\end{thm}
\begin{prf}
Since, by item $(iii)$ of \Cref{same as NV}, we have $\varepsilon_{\mathcal{Q}}^t(w) = \tau^{-1} E_2(we_2)$, it follows that $\varepsilon_{\mathcal{Q}}^t(vw)
= v\,\varepsilon_{\mathcal{Q}}^t(w)$ for all $v \in A'\cap A_1.$ Moreover, by \Cref{item (2)Delta(1)}, we know that $\Delta_{\mathcal{Q}}(wv)
= \Delta_{\mathcal{Q}}(w)(v\otimes 1), \text{ for all } w \in A'\cap A_2 \text{ and } v \in A'\cap A_1.$ Consequently, it is straightforward to verify that $({\mathcal{Q}}, \tilde{\Delta}_{\mathcal{Q}}, \tilde{\varepsilon}_{\mathcal{Q}})$ is a coalgebra. From \Cref{delta homo}, item $(i)$, it follows immediately that 
$\tilde{\Delta}_{\mathcal{Q}}$ is an algebra homomorphism. 
Furthermore, from \Cref{dfnnewimp} and \Cref{imp for kac alg 2}, it follows readily that for all 
$w, w', w'' \in A'\cap A_2$,
\[
\begin{aligned}
\tilde{\varepsilon}_{\mathcal{Q}}(w w' w'') 
= \tilde{\varepsilon}_{\mathcal{Q}}\!\big(w\, w'_{(\tilde{1})}\big)\, 
   \tilde{\varepsilon}_{\mathcal{Q}}\!\big(w'_{(\tilde{2})} w''\big) 
= \tilde{\varepsilon}_{\mathcal{Q}}\!\big(w\, w'_{(\tilde{2})}\big)\, 
   \tilde{\varepsilon}_{\mathcal{Q}}\!\big(w'_{(\tilde{1})} w''\big).
\end{aligned}
\]
In addition from \Cref{imp for kac algebra}, we have $$(\mathrm{id}\otimes \tilde{\Delta}_{\mathcal{Q}})\!\circ\!
\tilde{\Delta}_{\mathcal{Q}}(1)
= \big(\tilde{\Delta}_{\mathcal{Q}}(1)\otimes 1\big)\big(1\otimes 
\tilde{\Delta}_{\mathcal{Q}}(1)\big)
= \big(1\otimes \tilde{\Delta}_{\mathcal{Q}}(1)\big)
\big(\tilde{\Delta}_{\mathcal{Q}}(1)\otimes 1\big).
$$
Hence, $({\mathcal{Q}}, \tilde{\Delta}_{\mathcal{Q}}, \tilde{\varepsilon}_{\mathcal{Q}}, \tilde{S}_
{\mathcal{Q}})$
is a weak bialgebra. Using \Cref{delta homo}, item $(ii)$, it follows that for any 
$w \in {\mathcal{Q}}$,
\[
w_{(\tilde{1})}\, \tilde{S}_{\mathcal{Q}}(w_{(\tilde{2})})
= \tilde{\varepsilon}_{\mathcal{Q}}(1_{(\tilde{1})} w)\, 1_{(\tilde{2})}.
\]
By definition, $\tilde{S}_{\mathcal{Q}}$ is bijective, and it is 
clearly both an anti-algebra and anti-coalgebra map. Moreover, by \Cref{revision33}, $
\tilde{\varepsilon}_{\mathcal{Q}}\!\big(
\tilde{S}_{\mathcal{Q}}(w)\big)
= \tilde{\varepsilon}_{\mathcal{Q}}(w).$ Thus, for any $w \in {\mathcal{Q}}$,
\[
\begin{aligned}
\tilde{S}_{\mathcal{Q}}(w_{(\tilde{1})})\, w_{(\tilde{2})} 
&= \tilde{S}_{\mathcal{Q}}\Big(
   \tilde{\varepsilon}_{\mathcal{Q}}\big(1_{(\tilde{1})} 
   \tilde{S}_{\mathcal{Q}}^{-1}(w)\big)\, 1_{(\tilde{2})}\Big) \\[4pt]
&= \tilde{\varepsilon}_{\mathcal{Q}}\big(w\,
   \tilde{S}_{\mathcal{Q}}(1_{(\tilde{1})})\big)\, 
   \tilde{S}_{\mathcal{Q}}(1_{(\tilde{2})}) \\[4pt]
&= d^{-1}\,\sum_i \varepsilon_{\mathcal{Q}}(kw u_i^*)\, 
   r^-_1(u_i k^{-1}).
\end{aligned}
\]
Using the definition of a quasi-basis and applying the same line of argument, as it was in \Cref{Weak Kac}, it is straightforward to verify that
\[
d^{-1}\,\sum_i \varepsilon_{\mathcal{Q}}(kw u_i^*)\, 
r^-_1(u_i k^{-1})
= d^{-1}\,\sum_i \varepsilon_{\mathcal{Q}}(kw u_i)\, 
r^-_1(u_i^* k^{-1}).
\]
Thus, we obtain
\[
\begin{aligned}
\tilde{S}_{\mathcal{Q}}(w_{(\tilde{1})})\, w_{(\tilde{2})} 
&= \sum_id^{-1}\, \varepsilon_{\mathcal{Q}}(kw u_i)\, 
   r^-_1(u_i^* k^{-1}) \\[4pt]
&= \sum_i(\varepsilon_{\mathcal{Q}} \otimes \mathrm{id})\big(d^{-1}kwu_i \otimes 
   r^-_1(u_i^* k^{-1})\big).
\end{aligned}
\]
Now, observe that
\[
\begin{aligned}
\sum_i r^-_1(u_i^*)\, r^-_1(k^{-1}) \otimes d^{-1}kwu_i
&= (1 \otimes kw)\,\Delta_{\mathcal{Q}}(1)\,(r^-_1(k^{-1}) \otimes 1) \\
&= (1 \otimes kw)\,\Delta_{\mathcal{Q}}(1)\,(1 \otimes k^{-1}) \quad \text{ (by \Cref{imp form Dleta homo} item $(iv)$)}\\
&= \sum_i r^-_1(u_i^*) \otimes d^{-1}kwu_i k^{-1} \\
&=\sum_i r^-_1(u_i^*) \otimes d^{-1}kwk^{-1}u_i.
\end{aligned}
\]
Therefore,
\[
\begin{aligned}
\tilde{S}_{\mathcal{Q}}(w_{(\tilde{1})})\, w_{(\tilde{2})} 
&= \sum_i(\varepsilon_{\mathcal{Q}} \otimes \mathrm{id})\big(d^{-1}kwk^{-1}u_i \otimes 
   r^-_1(u_i^*)\big) \\[4pt]
&= \sum_i r^-_1(u_i^*)\, \varepsilon_{\mathcal{Q}}(d^{-1}kwk^{-1}u_i) \\[4pt]
&= 1_{(\tilde{1})}\, 
   \tilde{\varepsilon}_{\mathcal{Q}}(w\, 1_{(\tilde{2})}).
\end{aligned}
\]
Furthermore, as the map $\tilde{S}_{\mathcal{Q}}$ is both an 
anti-algebra and an anti-coalgebra map, \cite[Lemma~7.4]{Ni} yields
\[
\tilde{S}_{\mathcal{Q}}(w_{(\tilde{1})})\, 
w_{(\tilde{2})}\, 
\tilde{S}_{\mathcal{Q}}(w_{(\tilde{3})})
= \tilde{S}_{\mathcal{Q}}(w).
\]
Consequently, $({\mathcal{Q}}, \tilde{\Delta}_{\mathcal{Q}}, \tilde{\varepsilon}_{\mathcal{Q}}, \tilde{S}_{\mathcal{Q}})$ forms a weak Hopf algebra. Finally, using \Cref{item (2)Delta(1)} and \Cref{imp form Dleta homo}, item $(iv)$, it follows that $\tilde{\Delta}_{\mathcal{Q}}$ is $\dagger$-preserving. Hence, the proof is complete.  

\qed

\end{prf}

\begin{rmrk}\label{impfor mainth}
By \Cref{same as NV}, item ($iii$), it is straightforward to verify that $\varepsilon^t_{\mathcal{Q}}(k)=k.$
Hence, using \cite[Lemma 2.5]{BNS}, for any \(w\in A'\cap A_2\) we obtain
\begin{eqnarray}\nonumber
\tilde{\varepsilon}^t_{\mathcal{Q}}(w)
&=&\tilde{\varepsilon}_{\mathcal{Q}}(1_{(\tilde{1})}w)1_{(\tilde{2})}\\\nonumber
&=&{\varepsilon}_{\mathcal{Q}}(k1_{({1})}w)k^{-1}1_{({2})}\\\nonumber
&=&k^{-1}{\varepsilon}_{\mathcal{Q}}(1_{({1})}kw)1_{({2})}\\\nonumber
&=&k^{-1}{\varepsilon}^t_{\mathcal{Q}}(kw)\\\nonumber
&=&{\varepsilon}^t_{\mathcal{Q}}(w).
\end{eqnarray}
Therefore, $\tilde{\varepsilon}^t_{\mathcal{Q}}
=
\varepsilon^t_{\mathcal{Q}}.$ Furthermore, by \Cref{same as NV}, item ($iii$), we have $\mathcal{Q}_t
=
\tilde{\varepsilon}^t_{\mathcal{Q}}(\mathcal{Q})
=
A'\cap A_1.$ Moreover, by \Cref{property of reflection} and \Cref{dfnnewimp}, we obtain
\(\mathcal{Q}_s=\tilde{\varepsilon}^s_{\mathcal Q}(\mathcal Q)
=\tilde{S}_{\mathcal Q}(A'\cap A_1)
=A_1'\cap A_2\). It is also straightforward to verify that
\(\varepsilon^s_{\mathcal Q}(w)=\tilde{\varepsilon}^s_{\mathcal Q}(k^{-1}wk)\)
for every \(w\in A'\cap A_2\). Consequently,
\(\varepsilon^s_{\mathcal Q}(\mathcal Q)
=\tilde{\varepsilon}^s_{\mathcal Q}(\mathcal Q)
=A_1'\cap A_2\). 
\end{rmrk}

Using the the non-degenerate duality \(\langle \, \cdot , \cdot \, \rangle_{\mathrm{NV}}\) from \Cref{MAINTHM1}, we obtain the following result.
\begin{thm}\label{iiimpfor mainth}
Let $B\subset A$ be an inclusion of simple unital $C^*$-algebras with a conditional expectation of index-finite type. Furthermore, if the inclusion is of depth $2$, then ${\mathcal{P}}$ admits a weak $C^*$-Hopf algebra structure.
\end{thm}

An interesting feature of this weak \(C^*\)-Hopf algebra $\mathcal{P}$ is that its underlying \(C^*\)-algebra structure need not coincide with the usual \(C^*\)-algebra structure of the relative commutant \(B'\cap A_1\); see \Cref{discussion about difference} for further details.

\section{Action and crossed product}\label{Action and crossed product}
In this section, we study an inclusion $B \subset A$ of simple unital $C^*$-algebras with a conditional expectation of index-finite type, with the aim of characterizing when this inclusion has depth $2$. To this end, we first define a (left) action of ${\mathcal{Q}}$ on $A_1$. It follows that $A$ is the fixed-point subalgebra, while $A_2$ is isomorphic to the crossed product $A_1 \rtimes {\mathcal{Q}}.$ In the setting of simple $C^*$-algebras, the downward basic construction does not exist in general (see \cite[Remark 4.3]{I} and \cite{KajiwaraWatatani}). Hence, in contrast to the $II_1$ factor case, the above observation does not provide a direct characterization of the depth $2$ inclusion $B\subset A$. Nevertheless, it allows us to deduce the existence of a (left) action of ${\mathcal{P}}$ on $A$ such that $B$ is the corresponding fixed point subalgebra and $
A_1 \cong A \rtimes {\mathcal{P}}.$ First, we record the following propositions, whose proofs are analogous to those of Propositions~6.1 and~6.2 of \cite{NV}, and hence omitted.

\begin{ppsn}\label{action on V} Define $$\triangleright :
{\mathcal{Q}} \otimes A_1 \;\to\; A_1 \quad \text{ by } \quad w\triangleright x= \tau^{-1}E_2(wxe_2).$$
Then $\triangleright$ defines a (left) action of ${\mathcal{Q}}$ on $A_1$.
\end{ppsn}
\begin{ppsn} $A$ is the fixed point subalgebra of $A_1$, i.e., $
A^{{\mathcal{Q}}}_1 = A.$
\end{ppsn}
\begin{thm}
There exists a $*$-algebra isomorphism $\Theta: A_1 \rtimes {\mathcal{Q}} \rightarrow A_2$ satisfying $$\Theta([x\otimes w])= x{\tilde{S}_{\mathcal{Q}}(k)}^{\frac{1}{2}}\,w \,{\tilde{S}_{\mathcal{Q}}(k)}^{-\frac{1}{2}},$$
 for all $x\in A_1$ and $w\in {\mathcal{Q}}$.
\end{thm}
\begin{prf}
Using item $(i)$ of \Cref{impforaction}, from \cite[Proposition 6.3]{NV}, it follows that $\Theta$ is a well defined, involution-preserving homomorphism. In the setting of simple unital $C^*$-algebras, we next show that $\Theta$ is an isomorphism. To this end, we define a map $\Phi: A_2 \rightarrow A_1 \rtimes {\mathcal{Q}}$ by $$\Phi(w)=\sum_j \big[E_2\big(w (\mu^{\prime}_j)^*\big) \otimes {\tilde{S}_{\mathcal{Q}}(k)}^{-\frac{1}{2}}\, \mu^*_j\, {\tilde{S}_{\mathcal{Q}}(k)}^{\frac{1}{2}}\big] \quad \text{for each } w \in A_2,$$ where $\{ \mu_j , \mu^{\prime}_j\}$ are as in \Cref{AB=C}. Now we will show that $\Phi$ is the inverse map of $\Theta$. For all $x\in A_1$ and $w\in {\mathcal{Q}}$, we have the following:
\begin{eqnarray}\nonumber
\Phi \circ \Theta ([x \otimes w])&=& \Phi\big(x{\tilde{S}_{\mathcal{Q}}(k)}^{\frac{1}{2}}\,w \,{\tilde{S}_{\mathcal{Q}}(k)}^{-\frac{1}{2}}\big)\\\nonumber
&=&\sum_j \big[ E_2\big(x {\tilde{S}_{\mathcal{Q}}(k)}^{\frac{1}{2}}\,w \,{\tilde{S}_{\mathcal{Q}}(k)}^{-\frac{1}{2}} (\mu^{\prime}_j)^*\big) \otimes {\tilde{S}_{\mathcal{Q}}(k)}^{-\frac{1}{2}}\, \mu^*_j \,{\tilde{S}_{\mathcal{Q}}(k)}^{\frac{1}{2}}\big]\\\nonumber
&=& \sum_j\big[  xE_2\big({\tilde{S}_{\mathcal{Q}}(k)}^{\frac{1}{2}}\,w\, {\tilde{S}_{\mathcal{Q}}(k)}^{-\frac{1}{2}} \,(\mu^{\prime}_j)^*\big) \otimes {\tilde{S}_{\mathcal{Q}}(k)}^{-\frac{1}{2}}\, \mu^*_j\, {\tilde{S}_{\mathcal{Q}}(k)}^{\frac{1}{2}}\big]
\\\nonumber
&=&\sum_j \big[ x\big(E_2\big({\tilde{S}_{\mathcal{Q}}(k)}^{\frac{1}{2}}\,w\, {\tilde{S}_{\mathcal{Q}}(k)}^{-\frac{1}{2}} \,(\mu^{\prime}_j)^*\big)\triangleright 1\big) \otimes {\tilde{S}_{\mathcal{Q}}(k)}^{-\frac{1}{2}}\, \mu^*_j \,{\tilde{S}_{\mathcal{Q}}(k)}^{\frac{1}{2}}\big]
\\\nonumber
&=& \big[ x \otimes \sum_j E_2\big({\tilde{S}_{\mathcal{Q}}(k)}^{\frac{1}{2}}\,w\, {\tilde{S}_{\mathcal{Q}}(k)}^{-\frac{1}{2}}\, (\mu^{\prime}_j)^*\big)\,{\tilde{S}_{\mathcal{Q}}(k)}^{-\frac{1}{2}}\, \mu^*_j \,{\tilde{S}_{\mathcal{Q}}(k)}^{\frac{1}{2}}\big].
\end{eqnarray}
Now it is easy to observe the following:
\begin{eqnarray}\nonumber
&& \sum_jE_2\big({\tilde{S}_{\mathcal{Q}}(k)}^{\frac{1}{2}}\,w \,{\tilde{S}_{\mathcal{Q}}(k)}^{-\frac{1}{2}} (\mu^{\prime}_j)^*\big)\,{\tilde{S}_{\mathcal{Q}}(k)}^{-\frac{1}{2}}\, \mu^*_j \,{\tilde{S}_{\mathcal{Q}}(k)}^{\frac{1}{2}}\\\nonumber
&=& \sum_j{\tilde{S}_{\mathcal{Q}}(k)}^{-\frac{1}{2}}\, E_2\big({\tilde{S}_{\mathcal{Q}}(k)}^{\frac{1}{2}}\,w\, {\tilde{S}_{\mathcal{Q}}(k)}^{-\frac{1}{2}}\, (\mu^{\prime}_j)^*\big)\mu^*_j\, {\tilde{S}_{\mathcal{Q}}(k)}^{\frac{1}{2}}\\\nonumber
&=&{\tilde{S}_{\mathcal{Q}}(k)}^{-\frac{1}{2}}\, {\tilde{S}_{\mathcal{Q}}(k)}^{\frac{1}{2}}\,w\, {\tilde{S}_{\mathcal{Q}}(k)}^{-\frac{1}{2}}\, {\tilde{S}_{\mathcal{Q}}(k)}^{\frac{1}{2}}\\\nonumber
&=&w.
\end{eqnarray}
Hence we get $\Phi \circ \Theta ([x \otimes w])=[x \otimes w].$ Also for any $w\in {\mathcal{Q}}$, we obtain the following:
\begin{eqnarray}\nonumber
\Theta \circ \Phi (w)&=&\sum_j \Theta\big(\big[E_2\big(w \,(\mu^{\prime}_j)^*\big) \otimes {\tilde{S}_{\mathcal{Q}}(k)}^{-\frac{1}{2}}\, \mu^*_j\, {\tilde{S}_{\mathcal{Q}}(k)}^{\frac{1}{2}}\big]\big)\\\nonumber
&=&  \sum_jE_2\big(w \,(\mu^{\prime}_j)^*\big){\tilde{S}_{\mathcal{Q}}(k)}^{\frac{1}{2}}\,{\tilde{S}_{\mathcal{Q}}(k)}^{-\frac{1}{2}}\, \mu^*_j\, {\tilde{S}_{\mathcal{Q}}(k)}^{\frac{1}{2}}\, {\tilde{S}_{\mathcal{Q}}(k)}^{-\frac{1}{2}}\\\nonumber
&=& w.
\end{eqnarray}
Hence we are done.
\qed
\end{prf}
\medskip

It is easy to see that the action defined in \Cref{action on V} is minimal. From \Cref{iiimpfor mainth}, we know that the non-degenerate duality $\langle\,\cdot,\cdot\,\rangle_{\mathrm{NV}}$ endows \(\mathcal{P}\) with the structure of a weak \(C^*\)-Hopf algebra. We denote this weak \(C^*\)-Hopf algebra by $(\mathcal{P},\tilde{\Delta}_{\mathcal{P}},\tilde{\varepsilon}_{\mathcal{P}},\tilde{S}_{\mathcal{P}}).$ Furthermore, we denote by \(\tilde{\cdot}\) the multiplication on \(\mathcal{P}\) induced by duality from the weak \(C^*\)-Hopf algebra $(\mathcal{Q},\tilde{\Delta}_{\mathcal{Q}},\tilde{\varepsilon}_{\mathcal{Q}},\tilde{S}_{\mathcal{Q}}),$ and by \(\tilde{1}\) its unit element. 

We now define a left action of \(\mathcal{P}\) on \(A\) whose fixed point subalgebra is \(B\) and for which
\(A_1\cong A\rtimes\mathcal{P}\). To this end, we first establish several auxiliary results that will be used throughout the remainder of this section. Set $u=r_1^+(k^{-1})$. Then \(u\) is a positive invertible central element of \(B'\cap A\). Since \(r_1^+\) is an anti-isomorphism (see \Cref{aboutreflectionoperator}), it follows immediately that $r_1^+(k)=u^{-1}.$

\begin{lmma}\label{newproduct}
For all \(x,y\in\mathcal{P}\), $x\,\tilde{\cdot}\,y=d\,x\,u\,y.$ 
\end{lmma}
\begin{prf}
From the weak \(C^*\)-Hopf algebra structure on \(\mathcal{P}\) induced by the duality $\langle\,\cdot,\cdot\,\rangle_{\mathrm{NV}},$
we have
\begin{eqnarray}\nonumber
\langle\,x\,\tilde{\cdot}\,y, w\,\rangle_{\mathrm{NV}}&=&\langle\,x , w_{(\tilde{1})}\,\rangle_{\mathrm{NV}}\langle\,y , w_{(\tilde{2})}\,\rangle_{\mathrm{NV}}\\\label{needed1}
&=&d \,\langle\,x , w_{(1)}\,\rangle_{\mathrm{NV}}\langle\,y , k^{-1}w_{(2)}\,\rangle_{\mathrm{NV}}.
\end{eqnarray}
Furthermore,
\begin{eqnarray}\nonumber
\langle\,y , k^{-1}w_{(2)}\,\rangle_{\mathrm{NV}}&=&d\tau^{-2}\tr(ye_2e_1k^{-1}w_{(2)})\\\nonumber
&=&d\tau^{-2}\tr(ye_2e_1uw_{(2)}) \quad (\text{by \Cref{revision2}})\\\nonumber
&=&d\tau^{-2}\tr(ye_2e_1w_{(2)}u)\\\nonumber
&=&d\tau^{-2}\tr(uye_2e_1w_{(2)})\\\nonumber
&=& \langle\,uy , w_{(2)}\,\rangle_{\mathrm{NV}}
\end{eqnarray}
Thus from \Cref{needed1}, we have $
\langle\,x\,\tilde{\cdot}\,y, w\,\rangle_{\mathrm{NV}}
=d \,\langle\,x , w_{(1)}\,\rangle_{\mathrm{NV}}\langle\,uy , w_{(2)}\,\rangle_{\mathrm{NV}}=\langle\,d\,x\,u\,y, w\,\rangle_{\mathrm{NV}}.$
\qed
\end{prf}
\begin{lmma}\label{newunit}
The unit of $(\mathcal{P},\tilde{\Delta}_{\mathcal{P}},\tilde{\varepsilon}_{\mathcal{P}},\tilde{S}_{\mathcal{P}})$ is $\tilde{1}=\frac{u^{-1}}{d}.$ 
\end{lmma}
\begin{prf}From the duality, we have the following:
\begin{eqnarray}\nonumber
\langle\,\tilde{1}, w\,\rangle_{\mathrm{NV}}&=& \tilde{\varepsilon}_{\mathcal{Q}}(w)\\\nonumber
&=&d^{-1}\langle\,1, kw\,\rangle_{\mathrm{NV}}\\\nonumber
&=&\tau^{-2}\tr(e_2e_1u^{-1}w) \quad (\text{by \Cref{revision2}})\\\nonumber
&=&\tau^{-2}\tr(u^{-1} e_2e_1w) \\\nonumber
&=& \langle\,\frac{u^{-1}}{d}, w\,\rangle_{\mathrm{NV}}.
\end{eqnarray}
\qed
\end{prf}
\begin{rmrk}\label{discussion about difference}
By \Cref{newproduct} and \Cref{newunit}, the multiplication and unit on \(\mathcal{P}\) induced by the duality coincide with the usual multiplication and unit on the relative commutant \(B'\cap A_1\) if and only if \(k=\mathrm{Ind}_{W}(E_1|_{A'\cap A_1})\) is a scalar.
\end{rmrk}
\begin{lmma}\label{newantipode}
For every \(x\in\mathcal{P}\), $\tilde{S}_{\mathcal{P}}(x)=u^{-1}k^{-1}r_1^+(x).$ 
\end{lmma}
\begin{prf}
We have the following:
\begin{eqnarray}\nonumber
\langle\,\tilde{S}_{\mathcal{P}}(x), w\,\rangle_{\mathrm{NV}}&=&\langle\,x, \tilde{S}_Q(w)\,\rangle_{\mathrm{NV}}\\\nonumber
&=&\langle\,x, r_1^-(kwk^{-1})\,\rangle_{\mathrm{NV}}\\\nonumber
&=&\langle\,r_1^+(x), kwk^{-1}\,\rangle_{\mathrm{NV}}\quad (\text{by \Cref{impr-r+}})\\\nonumber
&=&d\tau^{-2}\tr(k^{-1}r_1^+(x) e_2e_1u^{-1}w) \quad (\text{by \Cref{revision2}})\\\nonumber
&=&d\tau^{-2}\tr(u^{-1}k^{-1}r_1^+(x) e_2e_1w) \\\nonumber
&=&\langle\,u^{-1}k^{-1}r_1^+(x), w\,\rangle_{\mathrm{NV}}.
\end{eqnarray}
\qed
\end{prf}

For every \(x\in\mathcal{P}\), write
\(\tilde{\Delta}_{\mathcal{P}}(x)=x_{(1)}\otimes x_{(2)}\).
We next establish the following lemma.
\begin{lmma}\label{revision3}
For every \(u'\in B'\cap A\) and every \(x\in\mathcal{P}\), $\tilde{\Delta}_{\mathcal{P}}(u'x)
=u'x_{(1)}\otimes x_{(2)}.$
\end{lmma}

\begin{prf}
We have the following:
\begin{eqnarray}\nonumber
\langle\,\tilde{\Delta}_{\mathcal{P}}(u'x), w\otimes w'\,\rangle_{\mathrm{NV}}&=&\langle\,u'x, ww'\,\rangle_{\mathrm{NV}}\\\nonumber
&=& d\tau^{-2}\tr(u'x e_2e_1ww')\\\nonumber
&=& d\tau^{-2}\tr(x e_2e_1u'ww')\\\nonumber
&=&\langle\,x, r^+_1(u')w w'\,\rangle_{\mathrm{NV}}\\\nonumber
&=&\langle\,x_{(1)}, r^+_1(u')w\,\rangle_{\mathrm{NV}}\langle\,x_{(2)}, w'\,\rangle_{\mathrm{NV}}\\\nonumber
&=&\langle\,u'x_{(1)}, w\,\rangle_{\mathrm{NV}}\langle\,x_{(2)}, w'\,\rangle_{\mathrm{NV}}.
\end{eqnarray}
\qed
\end{prf}

It is straightforward to verify that the \( * \)-structure on \(\mathcal{P}\), arising via the non-degenerate duality $\langle\,\cdot,\cdot\,\rangle_{\mathrm{NV}}$
as the dual of the weak \(C^*\)-Hopf algebra $(\mathcal{Q},\tilde{\Delta}_{\mathcal{Q}},\tilde{\varepsilon}_{\mathcal{Q}},\tilde{S}_{\mathcal{Q}}),$ coincides with the usual \( * \)-structure on the relative commutant \(\mathcal{P}\). Moreover, by taking adjoints in \Cref{revision3}, we obtain
\begin{equation}\label{revision3adj}
\tilde{\Delta}_{\mathcal{P}}(xu')
=x_{(1)}u'\otimes x_{(2)},\quad
\text{ for every }\,u'\in B'\cap A, \text{ and }\ x\in\mathcal{P}.
\end{equation}
The following result will also be useful.
\begin{lmma}\label{revision1308}
We have $\mathcal{P}_t=B'\cap A$ and 
$\mathcal{P}_s=u^{-1}(A'\cap A_1).$
\end{lmma}
\begin{prf}
For any $u'\in B'\cap A$ we have 
\begin{eqnarray}\nonumber
\tilde{\varepsilon}^t_{\mathcal{P}}(u')&=&\tilde{\varepsilon}_{\mathcal{P}}(\tilde{1}_{(1)}\, \tilde{\cdot}\,u')\tilde{1}_{(2)}\\\nonumber
&=&\tilde{\varepsilon}_{\mathcal{P}}({1}_{(1)} u^{-1}uu'){1}_{(2)}\quad (\text{by \Cref{newunit} and \Cref{revision3adj}})\\\nonumber
&=&(\tilde{\varepsilon}_{\mathcal{P}}\otimes \mathrm{id})\circ\tilde{\Delta}_{\mathcal{P}} (u')\quad (\text{by \Cref{revision3adj}})\\\nonumber
&=&u'.
\end{eqnarray}
Hence, \(B'\cap A\subseteq\mathcal{P}_t\). For the reverse inclusion, observe that \(\mathcal{P}_t\) is isomorphic to \(\mathcal{Q}_s=A_1'\cap A_2\). Moreover, by \cite[Theorem~3.27]{BGPS}, \(B'\cap A\) is isomorphic to \(A_1'\cap A_2\). Therefore, \(\dim(B'\cap A)=\dim(\mathcal{P}_t)\), and hence \(B'\cap A=\mathcal{P}_t\). Furthermore using \Cref{newantipode} and \Cref{revision1} we get $\mathcal{P}_s=u^{-1}(A'\cap A_1).$
\qed
\end{prf}

By \Cref{impfor mainth} and \Cref{revision1308}, both \(\mathcal{P}\) and \(\mathcal{Q}\) are biconnected weak \(C^*\)-Hopf algebras. We now proceed to define a left action of \(\mathcal{P}\) on \(A\).
\begin{ppsn}\label{MAINTHM2}
The map $ \triangleright : {\mathcal{P}} \otimes A \rightarrow A$ given by $$x \triangleright a= d^2\,u^{\frac{1}{2}}\,x_{(1)}\,u^{\frac{1}{2}}\, a\, u^{\frac{1}{2}}\,{\tilde{S}}_{\mathcal{P}}(x_{(2)})\,u^{\frac{1}{2}},$$ defines a (left) action of ${\mathcal{P}}$ on $A$. Furthermore, we have
$
A^{\mathcal{P}} = B.$

\end{ppsn}
To prove \Cref{MAINTHM2}, we begin by establishing a few preliminary results that will play a crucial role in the proof.

\begin{lmma}\label{revision8}
For every \(x\in\mathcal{P}\) and \(a\in A\), define $x\,\tilde{\triangleright}\,a
=
x_{(1)}\,a\,\tilde{S}_{\mathcal{P}}(x_{(2)}).$ 
Then $x\,\tilde{\triangleright}\,a\in A.$
\end{lmma}
\begin{prf}
The standard argument is analogous to those of Propositions 4.5 in \cite{KN}. However, in the present setting, we sketch the proof for completeness. We have to show that, for any $x\in {\mathcal{P}}$ and $a\in A$ we have $x\, \tilde{\triangleright}\, a\in A$. For that we consider the coaction dual to the action defined in \Cref{action on V}. Let $\rho:A_1 \rightarrow A_1 \otimes {\mathcal{P}}$ given by $\rho(y)= y^{(0)} \otimes y^{(1)},$ denote the dual coaction. Then for any $w\in {\mathcal{Q}}$ and $y\in A_1$ we have, $w\triangleright y=y^{(0)}\,\big<y^{(1)}, w\big>_{\text{NV}} $. Also for any $w\in {\mathcal{Q}}$, $x\in {\mathcal{P}}$, and $a\in A$ the following holds $
w \triangleright (a x) \;=\; a \,\langle x_{(2)}, w \rangle_{\mathrm{NV}}\, x_{(1)}$, and thus it is easy to observe that for any $x\in {\mathcal{P}},$ $\rho(x)$ is the comultiplication, i.e.,
$\rho(x)=x^{(0)} \otimes x^{(1)}= x_{(1)} \otimes x_{(2)}$. Since $A$ is the fixed point subalgebra of the action defined in \Cref{action on V}, the coinvariant subalgebra of the dual coaction is also $A$. Now we have the following:
\begin{eqnarray}\nonumber
\rho(x \, \tilde{\triangleright}\, a)&=& x_{(1)}\,a^{(0)}\,{\tilde{S}}_{\mathcal{P}}(x_{(4)}) \otimes x_{(2)}\,\tilde{\cdot}\,\tilde{\varepsilon}^{t}_{\mathcal{P}}(a^{(1)})\,\tilde{\cdot}\,{\tilde{S}}_{\mathcal{P}}(x_{(3)}).
\end{eqnarray}
Now in this setting, it is straightforward to verify that for any $x,z\in {\mathcal{P}}$ the following holds,
$
\tilde{\varepsilon}^{t}_{\mathcal{P}}(x \,\tilde{\cdot}\,z) = x_{(1)} \,\tilde{\cdot}\, \tilde{\varepsilon}^{t}_{\mathcal{P}}(z) \,\tilde{\cdot}\, {\tilde{S}}_{\mathcal{P}}(x_{(2)}),
$  
hence we have  
\begin{eqnarray}\nonumber
\rho(x  \,\tilde{\triangleright} \, a) = x_{(1)}\, a^{(0)}\,{\tilde{S}}_{\mathcal{P}}(x_{(3)}) \otimes \tilde{\varepsilon}^{t}_{\mathcal{P}}\big(x_{(2)}\,\tilde{\cdot}\, a^{(1)}\big) 
= (x  \,\tilde{\triangleright} \, a)^{(0)} \otimes \tilde{\varepsilon}^{t}_{\mathcal{P}}\big((x  \,\tilde{\triangleright} \,a)^{(1)}\big).
\end{eqnarray}  
Hence, we obtain that $x \,\tilde{\triangleright} \,a \in A$.\qed
\end{prf}

\begin{lmma}\label{revisedaction}
For all $x,x'\in {\mathcal{P}}$ and $a,a' \in A$, the following relations hold:
\begin{enumerate}
\item[$(i)$] $(x\,\tilde{\cdot}\,x')\triangleright a= x\triangleright(x'\triangleright a).$
\item[$(ii)$] $\tilde{1}\triangleright a=a$.
\item[$(iii)$] $x\triangleright (aa')= (x_{(1)}\triangleright a)(x_{(2)}\triangleright a').$
\item[$(iv)$] $(x\triangleright a)^*= {\tilde{S}_{\mathcal{P}}(x)}^*\triangleright a^*.$
\item[$(v)$] $x\triangleright 1= \tilde{\varepsilon}^{t}_{\mathcal{P}}(x)\triangleright 1= d\tilde{\varepsilon}^{t}_{\mathcal{P}}(x)u$.
\end{enumerate}
\end{lmma}
\begin{prf}
\begin{enumerate}[leftmargin=*]
\item[($i$)] A direct computation shows that
\begin{eqnarray}\nonumber
x\triangleright(x'\triangleright a)&=&x\triangleright(d^2\,u^{\frac{1}{2}}\,x^{\prime}_{(1)}\,u^{\frac{1}{2}}\, a\, u^{\frac{1}{2}}\,{\tilde{S}}_{\mathcal{P}}(x^{\prime}_{(2)})\,u^{\frac{1}{2}})\\\nonumber
&=&d^4\,u^{\frac{1}{2}}\,x_{(1)}\,u\,x^{\prime}_{(1)}\,u^{\frac{1}{2}}\, a\,u^{\frac{1}{2}}\,{\tilde{S}}_{\mathcal{P}}(x^{\prime}_{(2)})\, u\,{\tilde{S}}_{\mathcal{P}}(x_{(2)})\,u^{\frac{1}{2}}\\\nonumber
&=&d^2\,u^{\frac{1}{2}}\,(x_{(1)}\,\tilde{\cdot}\,x^{\prime}_{(1)})\,u^{\frac{1}{2}}\, a\,u^{\frac{1}{2}}\,{\tilde{S}}_{\mathcal{P}}(x_{(2)}\,\tilde{\cdot}\,x^{\prime}_{(2)})\,u^{\frac{1}{2}}\\\nonumber
&=& (x\,\tilde{\cdot}\,x')\triangleright a.
\end{eqnarray}
\item[($ii$)]
We have the following:
\begin{eqnarray}\nonumber
\tilde{1}\triangleright a&=&d^2\,u^{\frac{1}{2}}\,\tilde{1}_{(1)}\,u^{\frac{1}{2}}\, a\, u^{\frac{1}{2}}\,{\tilde{S}}_{\mathcal{P}}(\tilde{1}_{(2)})\,u^{\frac{1}{2}}\\\label{needed2}
&=&d^2\,u^{-\frac{1}{2}}\,u\,\tilde{1}_{(1)}\,u^{\frac{1}{2}}\, a\, u^{-\frac{1}{2}}\,u\,{\tilde{S}}_{\mathcal{P}}(\tilde{1}_{(2)})\,u^{\frac{1}{2}}.
\end{eqnarray}
Since \(\tilde{\Delta}_{\mathcal{P}}(\tilde{1})\in \mathcal{P}_s\otimes\mathcal{P}_t\), we have $u\tilde{1}_{(1)}\in A'\cap A_1.$ Hence from \Cref{needed2} we have,
 \begin{eqnarray}\nonumber
\tilde{1}\triangleright a&=&d\,a\, u^{\frac{1}{2}}\,(\tilde{1}_{(1)}\,\tilde{\cdot}\,{\tilde{S}}_{\mathcal{P}}(\tilde{1}_{(2)}))\,u^{\frac{1}{2}}\\\nonumber
&=& \,a\, u^{\frac{1}{2}}\,u^{-1}\,u^{\frac{1}{2}}\\\nonumber
&=&a.
\end{eqnarray}
\item[($iii$)]Observe that
\begin{eqnarray}\nonumber
(x_{(1)}\triangleright a)(x_{(2)}\triangleright a')&=&d^4\,u^{\frac{1}{2}}\,x_{(1)}\,u^{\frac{1}{2}}\, a\, u^{\frac{1}{2}}\,{\tilde{S}}_{\mathcal{P}}(x_{(2)})\,u\,x_{(3)}\,u^{\frac{1}{2}}\, a'\, u^{\frac{1}{2}}\,{\tilde{S}}_{\mathcal{P}}(x_{(4)})\,u^{\frac{1}{2}}\\\label{needed3}
&=&d^3\,u^{\frac{1}{2}}\,x_{(1)}\,u^{\frac{1}{2}}\, a\,u^{\frac{1}{2}}\,({\tilde{S}}_{\mathcal{P}}(x_{(2)})\,\tilde{\cdot}\,x_{(3)})\,u\,u^{-\frac{1}{2}}\, a'\, u^{\frac{1}{2}}\,{\tilde{S}}_{\mathcal{P}}(x_{(4)})\,u^{\frac{1}{2}}.
\end{eqnarray}
Since \({\tilde{S}}_{\mathcal{P}}(x_{(2)})\,\tilde{\cdot}\,x_{(3)}\in \mathcal{P}_s\), we have $({\tilde{S}}_{\mathcal{P}}(x_{(2)})\,\tilde{\cdot}\,x_{(3)})u\in A'\cap A_1.$ Hence from \Cref{needed3} we obtain,
\begin{eqnarray}\nonumber
(x_{(1)}\triangleright a)(x_{(2)}\triangleright a')
&=&d^3\,u^{\frac{1}{2}}\,x_{(1)}\,u^{\frac{1}{2}}\, aa'\, u^{\frac{1}{2}}\,({\tilde{S}}_{\mathcal{P}}(x_{(2)})\,\tilde{\cdot}\,x_{(3)})\,u\,{\tilde{S}}_{\mathcal{P}}(x_{(4)})\,u^{\frac{1}{2}}\\\nonumber
&=&d^2\,u^{\frac{1}{2}}\,x_{(1)}\,u^{\frac{1}{2}}\, aa'\, u^{\frac{1}{2}}\,{\tilde{S}}_{\mathcal{P}}(x_{(2)})\,u^{\frac{1}{2}}\\\nonumber
&=& x\triangleright (aa').
\end{eqnarray}
\item[($iv$)]Indeed, we have
\begin{eqnarray}\nonumber
{\tilde{S}_{\mathcal{P}}(x)}^*\triangleright a^*&=&d^2\,u^{\frac{1}{2}}\,{\tilde{S}_{\mathcal{P}}(x_{(2)})}^*\,u^{\frac{1}{2}}\, a^*\, u^{\frac{1}{2}}\,\tilde{S}_{\mathcal{P}}\big( {\tilde{S}_{\mathcal{P}}(x_{(1)})}^*\big)\,u^{\frac{1}{2}}\\\nonumber
&=&d^2\,u^{\frac{1}{2}}\,{\tilde{S}_{\mathcal{P}}(x_{(2)})}^*\,u^{\frac{1}{2}}\, a^*\, u^{\frac{1}{2}}\,\tilde{S}_{\mathcal{P}}\big( \tilde{S}_{\mathcal{P}}^{-1}(x^*_{(1)})\big)\,u^{\frac{1}{2}}\\\nonumber
&=&d^2\,u^{\frac{1}{2}}\,{\tilde{S}_{\mathcal{P}}(x_{(2)})}^*\,u^{\frac{1}{2}}\, a^*\, u^{\frac{1}{2}}\,x^*_{(1)}\,u^{\frac{1}{2}}\\\nonumber
&=&(x\triangleright a)^*.
\end{eqnarray}
\item[($v$)]We compute
\begin{eqnarray}\nonumber
x \triangleright 1&=&d^2\,u^{\frac{1}{2}}\,x_{(1)}\,u\,{\tilde{S}}_{\mathcal{P}}(x_{(2)})\,u^{\frac{1}{2}}\\\nonumber
&=& d \,u^{\frac{1}{2}}\,(x_{(1)}\,\tilde{\cdot}\,{\tilde{S}}_{\mathcal{P}}(x_{(2)}))\,u^{\frac{1}{2}}\\\label{needed4}
&=& d\tilde{\varepsilon}^{t}_{\mathcal{P}}(x)u.
\end{eqnarray}
The last equality follows since $x_{(1)}\,\tilde{\cdot}\,{\tilde{S}}_{\mathcal{P}}(x_{(2)})\in \mathcal{P}_t=B'\cap A$ and $u^{\frac{1}{2}}$ belongs to the center of $B'\cap A$. Furthermore, for every $x'\in \mathcal{P}_t$, we have \(\tilde{\Delta}_{\mathcal{P}}(x')=(x'\,\tilde{\cdot}\,\tilde{1}_{(1)}) \otimes \tilde{1}_{(2)}.\) Therefore,
\begin{eqnarray}\nonumber
\tilde{\varepsilon}^{t}_{\mathcal{P}}(x) \triangleright 1&=&d^2\,u^{\frac{1}{2}}\,(\tilde{\varepsilon}^{t}_{\mathcal{P}}(x)\, \tilde{\cdot}\, \tilde{1}_{(1)})\,u\,{\tilde{S}}_{\mathcal{P}}(\tilde{1}_{(2)})\,u^{\frac{1}{2}}\\\nonumber
&=&d^2\,u^{\frac{1}{2}}\,\tilde{\varepsilon}^{t}_{\mathcal{P}}(x)\, u\, (\tilde{1}_{(1)}\,\tilde{\cdot}\,{\tilde{S}}_{\mathcal{P}}(\tilde{1}_{(2)}))\,u^{\frac{1}{2}}\\\nonumber
&=&d\,u^{\frac{1}{2}}\,\tilde{\varepsilon}^{t}_{\mathcal{P}}(x)\, u\, u^{-1}\,u^{\frac{1}{2}}\\\label{needed5}
&=& d\tilde{\varepsilon}^{t}_{\mathcal{P}}(x)u.
\end{eqnarray}
The result now follows from \Cref{needed4,needed5}.\qed

\end{enumerate}
\end{prf}

\noindent\textbf{Proof of Proposition \ref{MAINTHM2}: }
We first prove that \(x\triangleright a\in A\) for every \(x\in\mathcal{P}\) and \(a\in A\). By \Cref{revision3},
$\tilde{\Delta}_{\mathcal{P}}(u^{\frac{1}{2}}x)
=u^{\frac{1}{2}}x_{(1)}\otimes x_{(2)}$. Therefore, by \Cref{revision8},
$(u^{\frac{1}{2}}x)\,\tilde{\triangleright}\,(u^{\frac{1}{2}} \,a\,u^{\frac{1}{2}})
=
u^{\frac{1}{2}}\,x_{(1)}\,u^{\frac{1}{2}}\,a\,u^{\frac{1}{2}}\,\tilde{S}_{\mathcal{P}}(x_{(2)})\in A$. Since \(u^{\frac12}\in A\), it follows that
$x\,\triangleright \,a= d^2\,u^{\frac{1}{2}}\,x_{(1)}\,u^{\frac{1}{2}}\, a\, u^{\frac{1}{2}}\,{\tilde{S}}_{\mathcal{P}}(x_{(2)})\,u^{\frac{1}{2}}\in A$. It now follows from \Cref{revisedaction} that the above formula defines a (left) action of the weak \(C^*\)-Hopf algebra \(\mathcal{P}\) on \(A\).

\medskip

We now show that $
A^{\mathcal{P}} = B.$ Let us assume that $b \in B$. Then, for any $x \in {\mathcal{P}}$, we have the following: 
\begin{eqnarray}\nonumber
x \triangleright b &=& d^2\,u^{\frac{1}{2}}\,x_{(1)}\,u^{\frac{1}{2}}\, b\, u^{\frac{1}{2}}\,{\tilde{S}}_{\mathcal{P}}(x_{(2)})\,u^{\frac{1}{2}}\\\nonumber
&=& d^2\,u^{\frac{1}{2}}\,x_{(1)}\,u\,{\tilde{S}}_{\mathcal{P}}(x_{(2)})\,u^{\frac{1}{2}}\,b\\\nonumber
&=& d\,\tilde{\varepsilon}^{t}_{\mathcal{P}} (x)\,u\,b.
\end{eqnarray}
The last equality follows since $\tilde{\varepsilon}^{t}_{\mathcal{P}} (x)\in \mathcal{P}_t=B'\cap A$ and $u^{\frac{1}{2}}$ belongs to the center of $B'\cap A$. Moreover, $\tilde{\varepsilon}^{t}_{\mathcal{P}} (x)\triangleright b=d\,\tilde{\varepsilon}^{t}_{\mathcal{P}} (x)\,u\,b$. Therefore, $B \subset A^{\mathcal{P}}$. Now we will show the other part of the inclusion, i.e., $A^{\mathcal{P}} \subset B$. To this end, we first verify that for any $a \in A^{\mathcal{P}}$ and $x \in {\mathcal{P}}$, the following holds:
\begin{eqnarray}\nonumber
u^{-\frac{1}{2}}\,\,a\,u^{\frac{1}{2}}\,{\tilde{S}_{\mathcal{P}}(x)}\,u^{\frac{1}{2}}&=& d^2u^{-\frac{1}{2}}  a\,u^{\frac{1}{2}}\,\tilde{S}_{\mathcal{P}}(x_{(1)})\,u\,x_{(2)}\,u\,\tilde{S}_{\mathcal{P}}(x_{(3)}) \,u^{\frac{1}{2}}\\\nonumber
&=& d^2\,\tilde{S}_{\mathcal{P}}(x_{(1)})\,u^{\frac{1}{2}}\,\,u^{\frac{1}{2}}\, x_{(2)}\,u^{\frac{1}{2}}\, a\,u^{\frac{1}{2}}\,\tilde{S}_{\mathcal{P}}(x_{(3)}) \,u^{\frac{1}{2}}\\\nonumber
&=& \tilde{S}_{\mathcal{P}}(x_{(1)})\,u^{\frac{1}{2}}\,(x_{(2)} \triangleright a)\\\nonumber
&=& \tilde{S}_{\mathcal{P}}(x_{(1)})\,u^{\frac{1}{2}}\,(\tilde{\varepsilon}^{t}_{\mathcal{P}} (x_{(2)}) \triangleright a)\\\label{needed9}
&=&d^3\, \tilde{S}_{\mathcal{P}}(x_{(1)})\,u\,\tilde{\varepsilon}^{t}_{\mathcal{P}} (x_{(2)})\,u\,\tilde{1}_{(1)}\,u^{\frac{1}{2}}\, a\, u^{\frac{1}{2}}\,{\tilde{S}}_{\mathcal{P}}(\tilde{1}_{(2)})\,u^{\frac{1}{2}}.
\end{eqnarray}
Since \(\tilde{\Delta}_{\mathcal{P}}(\tilde{1})\in \mathcal{P}_s\otimes\mathcal{P}_t\), we have $\tilde{1}_{(1)}u\in A'\cap A_1.$ Hence from \Cref{needed9} we have,
\begin{eqnarray}\nonumber
u^{-\frac{1}{2}}\,\,a\,u^{\frac{1}{2}}\,{\tilde{S}_{\mathcal{P}}(x)}\,u^{\frac{1}{2}}
&=&d^3\, \tilde{S}_{\mathcal{P}}(x_{(1)})\,u\,\tilde{\varepsilon}^{t}_{\mathcal{P}} (x_{(2)})\,u\,\tilde{1}_{(1)}\,u\,u^{-\frac{1}{2}}\, a\, u^{\frac{1}{2}}\,{\tilde{S}}_{\mathcal{P}}(\tilde{1}_{(2)})\,u^{\frac{1}{2}}\\\nonumber
&=&d^2\, \tilde{S}_{\mathcal{P}}(x_{(1)})\,u\,\tilde{\varepsilon}^{t}_{\mathcal{P}} (x_{(2)})\,u^{\frac{1}{2}}\, a\, u^{\frac{1}{2}}\,(\tilde{1}_{(1)}\,\tilde{\cdot}\,{\tilde{S}}_{\mathcal{P}}(\tilde{1}_{(2)}))\,u^{\frac{1}{2}}\\\nonumber
&=&\, (\tilde{S}_{\mathcal{P}}(x_{(1)})\,\tilde{\cdot}\,x_{(2)}\,\tilde{\cdot}\,\tilde{S}_{\mathcal{P}}(x_{(3)}))\,u^{\frac{1}{2}}\, a\, u^{\frac{1}{2}}\,u^{-1}\,u^{\frac{1}{2}}\\\nonumber
&=&\, \tilde{S}_{\mathcal{P}}(x)\,u^{\frac{1}{2}}\, a.
\end{eqnarray}
Hence,
$u^{\frac{1}{2}}\,\tilde{S}_{\mathcal{P}}(x)\,u^{\frac{1}{2}}\,a=a\,u^{\frac{1}{2}}\,\tilde{S}_{\mathcal{P}}(x)\,u^{\frac{1}{2}}.$
Moreover, setting \(x=\tilde{S}_{\mathcal P}^{-1}(e_1)\), we obtain, for every for every \(a\in A^{\mathcal{P}}\),
$u^{\frac{1}{2}}\,e_1\,u^{\frac{1}{2}}\,a=a\,u^{\frac{1}{2}}\,e_1\,u^{\frac{1}{2}}.$ Hence
$u^{\frac{1}{2}}\,e_1\,u^{\frac{1}{2}}\,a\,e_1=a\,u^{\frac{1}{2}}\,e_1\,u^{\frac{1}{2}}\,e_1
= a\,u^{\frac{1}{2}}\,E_0(u^{\frac{1}{2}})\,e_1$.
Thus,
$u^{\frac{1}{2}}\,E_0(u^{\frac{1}{2}}\,a)
= a\,u^{\frac{1}{2}}\,E_0(u^{\frac{1}{2}})$.
Since \(u^{\frac12}\) is positive and invertible, \(E_0(u^{\frac12})\) is also invertible. Moreover, as \(u^{\frac12}\in B'\cap A\),
\(a=E_0(u^{\frac12}a)E_0(u^{\frac12})^{-1}\).
Hence \(a\in B\), completing the proof.
\qed

\begin{thm}\label{MAINTHM3}
There exists a $*$-algebra isomorphism $\Psi: A \rtimes {\mathcal{P}} \rightarrow A_1$, satisfying $\Psi([a\otimes x])=d \,a\,u^{\frac{1}{2}}\,x\,u^{\frac{1}{2}}\,,$ for all $a \in A$ and $x\in {\mathcal{P}}$. 
\end{thm}
\begin{prf}
The proof is inspired by the proof of \cite[Theorem 4.6]{KN}. The map \(\Psi\) is indeed well-defined, as one can check that $\Psi\big([a(\tilde{\varepsilon}^t_{\mathcal P}(y)\triangleright 1)\otimes x]\big)
=
\Psi\big([a\otimes(\tilde{\varepsilon}^t_{\mathcal P}(y)\,\tilde{\cdot}\,x)]\big)$ for all \(a\in A\) and \(x,y\in\mathcal{P}\). Since $B \subset A$ is an inclusion of simple unital $C^*$-algebras with a conditional expectation of index-finite type, and the inclusion is of depth $2$, following the same argument as in \Cref{AB=C}, we can find a family $\{\gamma_j, \gamma^{\prime}_j\} \subset {\mathcal{P}}$ such that any $x \in A_1$ can be written as
$
x = \sum_jE_1\big(x (\gamma^{\prime}_j)^*\big) \gamma^*_j.
$
We now define a map
$$
\varphi: A_1 \rightarrow A \rtimes {\mathcal{P}}, \qquad 
\varphi(x) = d^{-1}\sum_j\big[E_1\big(x\,u^{-\frac{1}{2}}\, (\gamma^{\prime}_j)^*\big) \otimes \,u^{-\frac{1}{2}}\,\gamma^*_j\big], \quad \text{for each } x \in A_1.
$$
Next, we show that $\varphi$ is the inverse map of $\Psi$:
\begin{eqnarray}\nonumber
\varphi \circ \Psi([a \otimes x]) 
&=& d\,\varphi( \,a\,u^{\frac{1}{2}}\,x\,u^{\frac{1}{2}}) \\\nonumber
&=&\sum_j \big[E_1\big(a\,u^{\frac{1}{2}}\, x (\gamma^{\prime}_j)^*\big) \otimes \,u^{-\frac{1}{2}}\, \gamma^*_j\big] \\\nonumber
&=& \sum_j\big[a\,u^{\frac{1}{2}}\, E_1\big(x (\gamma^{\prime}_j)^*\big) \otimes \,u^{-\frac{1}{2}}\,\gamma^*_j\big]\\\nonumber
&=& \sum_j\big[a\,u^{-\frac{1}{2}}\, E_1\big(x (\gamma^{\prime}_j)^*\big)u \otimes \,u^{-\frac{1}{2}}\, \gamma^*_j\big].
\end{eqnarray}
Since ${\mathcal{P}}_t= B' \cap A$ and $u^{-\frac{1}{2}}\, E_1\big(x (\gamma^{\prime}_j)^*\big)\in B'\cap A$, from \Cref{revisedaction} item ($v$), it follows that
\begin{eqnarray}\nonumber
\varphi \circ \Psi([a \otimes x]) 
&=&d^{-1} \sum_j\big[a\big(\big(u^{-\frac{1}{2}}\, E_1\big(x (\gamma^{\prime}_j)^*\big)\big)\triangleright 1\big) \otimes \,u^{-\frac{1}{2}}\, \gamma^*_j\big]\\\nonumber
&=&d^{-1} \sum_j\big[a\otimes \big(u^{-\frac{1}{2}}\, E_1\big(x (\gamma^{\prime}_j)^*\big)\big)\,\tilde{\cdot}\,(u^{-\frac{1}{2}}\, \gamma^*_j)\big]\\\nonumber
&=& \sum_j\big[a\otimes E_1\big(x (\gamma^{\prime}_j)^*\big)\,\gamma^*_j\big]\\\nonumber
&=& [a \otimes x].
\end{eqnarray}
 Conversely, we have 
 $\Psi \circ \varphi(x) =\sum_jE_1\big(x\,u^{-\frac{1}{2}}\, (\gamma^{\prime}_j)^*\big) \gamma^*_j\,u^{\frac{1}{2}}=x$. Hence, $\Psi$ is a linear isomorphism. We next show that \(\Psi\) is an algebra homomorphism.
 
 \begin{eqnarray}\nonumber
 \Psi([a\otimes x][a'\otimes x'])&=&\Psi([\, a (x_{(1)} \triangleright a') \otimes x_{(2)}\,\tilde{\cdot}\, x' \,])\\\nonumber
 &=& d\,a\, (x_{(1)} \triangleright a')\,u^{\frac{1}{2}}\,(x_{(2)}\,\tilde{\cdot}\, x')\,u^{\frac{1}{2}}\\\label{needed11}
 &=& d^2\,a \,(x_{(1)} \triangleright a')\,u^{\frac{1}{2}}\,x_{(2)}\,u\, x'\,u^{\frac{1}{2}}
 \end{eqnarray}
 We now observe that
 \begin{eqnarray}\nonumber
\bigl(x_{(1)} \triangleright a' \bigr)\,u^{\frac{1}{2}}\, x_{(2)}\,u 
&=& d^2\,u^{\frac{1}{2}}\,x_{(1)}\,u^{\frac{1}{2}}\, a'\, u^{\frac{1}{2}}\,{\tilde{S}}_{\mathcal{P}}(x_{(2)})\,u\, x_{(3)}\,u  \\\nonumber
&=& d\,u^{\frac{1}{2}}\,x_{(1)}\,u^{\frac{1}{2}}\, a'\, u^{-\frac{1}{2}}\,u\,({\tilde{S}}_{\mathcal{P}}(x_{(2)})\,\tilde{\cdot}\, x_{(3)})\,u  \\\nonumber
&=& u^{\frac{1}{2}}\,(x_{(1)}\,\tilde{\cdot}\,{\tilde{S}}_{\mathcal{P}}(x_{(2)})\,\tilde{\cdot}\, x_{(3)})\,u^{\frac{1}{2}}\, a'\, u^{\frac{1}{2}} \\\label{needed12}
&=& u^{\frac{1}{2}}\,x\,u^{\frac{1}{2}}\, a'\, u^{\frac{1}{2}} \\\nonumber
\end{eqnarray}
 Combining \Cref{needed11} and \Cref{needed12}, we obtain
 \begin{eqnarray}\nonumber
 \Psi([a\otimes x][a'\otimes x'])
 &=& d^2\,a \,u^{\frac{1}{2}}\,x\,u^{\frac{1}{2}}\, a'\, u^{\frac{1}{2}}\, x'\,u^{\frac{1}{2}}= \Psi([a\otimes x])\Psi([a'\otimes x']).
 \end{eqnarray}
 It remains to show that \(\Psi\) is \( * \)-preserving. Indeed, we have
 \begin{eqnarray}\nonumber
 \Psi({[a \otimes x]}^*)
&=& \Psi([\, (x_{(1)}^* \triangleright a^*) \otimes x_{(2)}^* \,])\\\nonumber
&=& d\,(x_{(1)}^* \triangleright a^*)\,u^{\frac{1}{2}}\,x_{(2)}^*\,u^{\frac{1}{2}}\\\nonumber
&=& d\,u^{\frac{1}{2}}\,x^*\,u^{\frac{1}{2}}\, a^*\, u^{\frac{1}{2}}\,u^{-\frac{1}{2}}  \quad (\text{by \Cref{needed12}})\\\nonumber
&=& {\Psi([a \otimes x])}^*. 
 \end{eqnarray}
 This completes the proof.
 \qed
\end{prf}
\begin{rmrk}
The action described in \Cref{MAINTHM2} is minimal. If the inclusion $B \subset A$ is irreducible, then ${\mathcal{P}}$ is a finite-dimensional Kac algebra, yielding the result of \cite{I}.
\end{rmrk}

\begin{dfn}[\cite{BG}]
Let \( B \subset A \) be a unital inclusion of \( C^* \)-algebras.  
The unitary normalizer of \( B \) in \( A \) is defined by $\mathcal{N}_A(B) = \{\, u \in \mathcal{U}(A) : uBu^* = B \,\}.$ The inclusion \( B \subset A \) is called regular if $\mathcal{N}_A(B)$ generates the \( C^* \)-algebra $A$.
\end{dfn}
 It was conjectured in \cite{BG2} and later proved in \cite{CKP} that any finite index regular inclusion of $II_1$ factors has depth at most 2. In the $C^*$-algebraic setting, it was proved in~\cite[Theorem~3.11]{BG} that any regular inclusion \( B \subset A \) of simple unital \( C^* \)-algebras with a conditional expectation of index-finite type has depth at most~\(2\). As a consequence, we obtain the following corollary.

\begin{crlre}
Let $B \subset A$ be a regular inclusion of simple unital $C^*$-algebras with a conditional expectation of index-finite type. Then there exists an action of the weak $C^*$-Hopf algebra ${\mathcal{P}} $ on $A$ such that $B$ is the fixed-point subalgebra and $A_1 \cong A \rtimes {\mathcal{P}}$.
\end{crlre}

\section*{Acknowledgement}
The author is grateful to Prof.~Keshab Chandra Bakshi for valuable discussions and guidance. The author also thanks the anonymous referee for various helpful comments and suggestions, which significantly improved the exposition of this paper.

\bigskip

\noindent \\
       {\em Department of Mathematics and Statistics,\\
         Indian Institute of Technology Kanpur,\\
         Uttar Pradesh 208016, India.}\\
        { Email adress: \texttt{biplabpal32@gmail.com, bpal21@iitk.ac.in}}

\bigskip


\bigskip
	\begin{thebibliography}{ABC}

	
	
\bibitem{BGoP}
K.~C. Bakshi, D. Goswami and B. Pal, Weak quantum hypergroups from finite index $C^*$-inclusions, {\em Int. Math. Res. Not. IMRN} {\bf 2026}, no.~10, Paper No. rnag087, 24 pp.; MR5069553



\bibitem{BGP}
K.~C. Bakshi, S. Guin and B. Pal, Von Neumann entropy of the angle operator between a pair of intermediate subalgebras, {\em J. Noncommut. Geom.} {\bf 20} (2026), no.~2, 421--433; MR5047636

\bibitem{BGPS}
K.~C. Bakshi, S. Guin, B. Pal and Sruthymurali, Higher reflections and entropy of canonical shifts for inclusions of $C^{*}$-algebras with finite Watatani index, {\em Internat. J. Math.} {\bf 37} (2026), no.~1, Paper No. 2650007; MR5002824

\bibitem{BGS}
K.~C. Bakshi, S. Guin and Sruthymurali, Fourier-theoretic inequalities for inclusions of simple $C^*$-algebras, {\em New York J. Math.} {\bf 29} (2023), 335--362; MR4564011	
	
\bibitem{BakshiVedlattice}
K.~C. Bakshi and V.~P. Gupta, Lattice of intermediate subalgebras, {\em  J. Lond. Math. Soc.} (2) {\bf 104} (2021), no.~5, 2082--2127; MR4368671

\bibitem{BG2}
K.~C. Bakshi and V.~P. Gupta, A few remarks on Pimsner-Popa bases and regular subfactors of depth 2, {\em Glasg. Math. J.} {\bf 64} (2022), no.~3, 586--602; MR4462379

\bibitem{BG}
K.~C. Bakshi and V.~P. Gupta, Regular inclusions of simple unital $C^*$-algebras, {\em M\"unster J. Math.} {\bf 18} (2025), no.~1, 181--200; MR4976177



\bibitem{B}
D.~H. Bisch, Bimodules, Higher relative commutants and the fusion algebra associated to a subfactor, {\em Operator algebras and their applications (Waterloo, ON, 1994/1995)}, 13--63, Fields Inst. Commun., 13, Amer. Math. Soc., Providence, RI, ; MR1424954


\bibitem{BNS}
G. B\"ohm, F. Nill and K. Szlach\'anyi, Weak Hopf algebras. I. Integral theory and $C^*$-structure, {\em J. Algebra} {\bf 221} (1999), no.~2, 385--438; MR1726707

\bibitem{BS}
G. B\"ohm and K. Szlach\'anyi, A coassociative $C^*$-quantum group with nonintegral dimensions, {\em Lett. Math. Phys.} {\bf 38} (1996), no.~4, 437--456; MR1421688

\bibitem{CKP}
J. Crann, D.~W. Kribs and R.~J. Pereira, Orthogonal unitary bases and a subfactor conjecture, {\em Proc. Amer. Math. Soc.} {\bf 151} (2023), no.~9, 3793--3799; MR4607624

\bibitem{D}
M.-C. David, Paragroupe d'Adrian Ocneanu et alg\`ebre de Kac, {\em Pacific J. Math.} {\bf 172} (1996), no.~2, 331--363; MR1386622

\bibitem{ER}
S. Echterhoff and M. R\o rdam, Inclusions of $C^*$-algebras arising from fixed-point algebras, {\em Groups Geom. Dyn.} {\bf 18} (2024), no.~1, 127--145; MR4705629


\bibitem{EN}
M. Enock and R. Nest, Irreducible inclusions of factors, multiplicative unitaries, and Kac algebras, {J. Funct. Anal.} {\bf 137} (1996), no.~2, 466--543; MR1387518

\bibitem{EV}
M. Enock and J.-M. Vallin, Inclusions of von Neumann algebras, and quantum groupoids, {\em J. Funct. Anal.} {\bf 172} (2000), no.~2, 249--300; MR1753177


\bibitem{FI}
F. Fidaleo and T. Isola, The canonical endomorphism for infinite index inclusions, {\em Z. Anal. Anwendungen} {\bf 18} (1999), no.~1, 47--66; MR1681843

\bibitem{I}
M. Izumi, Inclusions of simple $C^\ast$-algebras, {\em J. Reine Angew. Math.} {\bf 547} (2002), 97--138; MR1900138

\bibitem{JOPT}
J.~A. Jeong, H. Osaka, N.~C. Phillips and T. Teruya, Cancellation for inclusions of $C^\ast$-algebras of finite depth, {\em Indiana Univ. Math. J.} {\bf 58} (2009), no.~4, 1537--1564; MR2542972

\bibitem{J}
V.~F.~R. Jones, Index for subfactors, {\em Invent. Math.} {\bf 72} (1983), no.~1, 1--25; MR0696688

\bibitem{KN}
L. Kadison and D.~A. Nikshych, Frobenius extensions and weak Hopf algebras, {\em J. Algebra} {\bf 244} (2001), no.~1, 312--342; MR1856540

\bibitem{KajiwaraWatatani}
T. Kajiwara and Y. Watatani, Jones index theory by Hilbert $C^*$-bimodules and $K$-theory, {\em Trans. Amer. Math. Soc.} {\bf 352} (2000), no.~8, 3429--3472; MR1624182

\bibitem{L}
R. Longo, A duality for Hopf algebras and for subfactors. I, {\em Comm. Math. Phys.} {\bf 159} (1994), no.~1, 133--150; MR1257245

\bibitem{M}
M. Mukohara, Inclusions of simple $\rm C^\ast$-algebras arising from compact group actions, {\em J. Funct. Anal.} {\bf 288} (2025), no.~2, Paper No. 110702, 33 pp.; MR4813137

\bibitem{N}
D.~A. Nikshych, Duality for actions of weak Kac algebras and crossed product inclusions of $\rm II_1$ factors, {\em J. Operator Theory} {\bf 46} (2001), no.~3, 635--655; MR1897159


\bibitem{NV}
D.~A. Nikshych and L.~I. Vainerman, A characterization of depth 2 subfactors of ${\rm II}_1$ factors, {\em J. Funct. Anal.} {\bf 171} (2000), no.~2, 278--307; MR1745634


\bibitem{Ni}
F. Nill, Axioms for weak bialgebras, math.QA/9805104 (1998).

\bibitem{NSW}
F. Nill, K. Szlachanyi, and H.W. Wiesbrock, Weak Hopf algebras and reducible Jones inclusions of depth 2. I: From crossed products to Jones towers. arXiv preprint math/9806130 (1998).

\bibitem{R}
M. R\o rdam, Irreducible inclusions of simple $C^*$-algebras, {\em Enseign. Math.} {\bf 69} (2023), no.~3-4, 275--314; MR4599249

\bibitem{O}
A. Ocneanu, Quantized groups, string algebras and Galois theory for algebras, {\em Operator algebras and applications}, Vol. 2, 119--172, London Math. Soc. Lecture Note Ser., 136, Cambridge Univ. Press, Cambridge, ; MR0996454

\bibitem{SPreprint}
K. Szlach\'anyi, Weak Hopf algebra symmetries of $C^*$-algebra inclusions, in {\it Geometry Seminars, 2000 (Italian) (Bologna, 1999/2000)}, 229--246, Univ. Stud. Bologna, Bologna, ; MR1878130


\bibitem{S}
W. Szyma\'nski, Finite index subfactors and Hopf algebra crossed products, {\em Proc. Amer. Math. Soc.} {\bf 120} (1994), no.~2, 519--528; MR1186139

\bibitem{Watataniindex}
Y. Watatani, Index for $C^*$-subalgebras, {\em Mem. Amer. Math. Soc.} {\bf 83} (1990), no.~424, {\rm vi}+117 pp.; MR0996807



\end{thebibliography}
\end{document}